\documentclass[11pt]{article}
\usepackage[utf8]{inputenc}
\usepackage[T1]{fontenc}
\usepackage[frenchb,english]{babel} 
\usepackage{textcomp}
\usepackage{amsmath,amssymb}
\usepackage{amsthm}
                  
\usepackage{lmodern}

\usepackage[a4paper,margin=2cm]{geometry}

\usepackage{graphicx}             
\usepackage{xcolor}               
\usepackage{microtype,stmaryrd}          

\usepackage{hyperref}
\hypersetup{pdfstartview=XYZ}

\usepackage{bbm}
\usepackage{mathrsfs}
\usepackage{subfig}

\def\build#1_#2^#3{\mathrel{
\mathop{\kern 0pt#1}\limits_{#2}^{#3}}}

\def\midd{\,|\,}

\newtheorem{theorem}{Theorem}
\newtheorem{proposition}[theorem]{Proposition}
\newtheorem{definition}[theorem]{Definition}
\newtheorem{lemma}[theorem]{Lemma}

\def\w{\mathrm{w}}
\def\t{\mathcal{T}}
\def\b{\mathcal{B}}

\def\g{\mathcal{G}}
\def\W{\mathcal{W}}
\def\S{\mathcal{S}}
\def\T{\mathbb{T}}
\def\N{\mathbb{N}}
\def\M{\mathbb{M}}

\def\P{\mathbb{P}}
\def\E{\mathbb{E}}

\def\R{\mathbb{R}}

\def\n{\mathcal{N}}

\def\ve{{\varepsilon}}
\def\la{\longrightarrow}

\def\ov{\overline}

\def\dd{\mathrm{d}}
\def\wh{\widehat}
\def\wt{\widetilde}

\def\bp{\mathbf{p}}

\def\v{\mathcal{V}}

\def\bm{\mathbf{m}}
\def\bx{\mathbf{x}}

\def\rem{\noindent{\bf Remark. }}

\author{Jean-Fran\c cois Le Gall}
\title{The volume measure of the Brownian sphere is a  Hausdorff measure\footnote{Supported by the ERC Advanced Grant 740943 {\sc GeoBrown}}}
\date{\small Universit\'e Paris-Saclay}

\begin{document}

\maketitle

\begin{abstract}
We prove that the volume measure of the Brownian sphere is
equal to a constant multiple of the Hausdorff measure associated with
the gauge function $h(r)=r^4\log\log(1/r)$. This shows in particular
that the volume measure of the Brownian sphere is determined 
by its metric structure. 
\end{abstract}

\section{Introduction}

This work is concerned with the Brownian sphere, also known as the Brownian map,
which is a random compact metric space providing a universal model
of two-dimensional random geometry. The Brownian sphere has been shown to
be the scaling limit in the Gromov-Hausdorff sense of many different classes of
random planar maps such as triangulations or quadrangulations of
the sphere (see in particular \cite{Abr,AA,AA2,BJM,Uniqueness,Mar,Mar2,Mie-Acta}).
The first construction the Brownian sphere \cite{MM05,Invent} relied
on Brownian motion indexed by the Brownian tree, or equivalently on
the Brownian snake driven by a normalized Brownian excursion. In this construction,
the Brownian sphere is constructed as a quotient space of the interval $[0,1]$ 
and is naturally equipped with a volume measure defined as the pushforward
of Lebesgue measure on $[0,1]$ under the canonical projection. This volume measure 
appears as the limit of (scaled) counting measures on vertices when the Brownian sphere 
is written as the Gromov-Hausdorff-Prokhorov limit of large random planar
maps (see \cite[Theorem 7]{Disks} for the case of quadrangulations and \cite[Theorem 1.2]{Mar2}
for the more general case of bipartite planar maps with a prescribed degree sequence).
We mention that a very different approach to the Brownian sphere, involving deep connections with
Liouville quantum gravity has been developed by Miller and Sheffield
in a series of papers \cite{MS1,MS2,MS3,MS4}. More recently, 
Ding, Dub\'edat, Dunlap and Falconet \cite{DDDF} have studied Liouville first-passage
percolation metrics associated with mollified versions of the Gaussian free field and
were able to prove the tightness of these renormalized metrics. Gwynne and Miller \cite{GM} later
proved the uniqueness of the limit, which is determined by the Gaussian free field
and called the Liouville quantum gravity metric. A particular case of this metric
should correspond to a variant of the Brownian sphere metric.

The Brownian sphere is known to be homeomorphic to the $2$-sphere, but
its Hausdorff dimension is equal to $4$ \cite{Invent}. The original motivation of the
present work was to determine an exact Hausdorff measure function 
for the Brownian sphere. Our main result solves this problem, and 
also provides a natural interpretation of the volume measure.

We denote the Brownian sphere by $\bm_\infty$, and write $\mathrm{Vol}$
for the volume measure on $\bm_\infty$. For any gauge function $h$,
 we denote the associated Hausdorff measure by $m^h$.

\begin{theorem}
\label{main}
For every $r\in(0,1/4)$, set $h(r)=r^4\log\log(1/r)$.
There exists a constant $\kappa>0$ such that we have almost surely for
every Borel subset $A$ of $\bm_\infty$,
$$m^h(A)= \kappa\,\mathrm{Vol}(A).$$
\end{theorem}

It seems hopeless to compute the exact value of the constant $\kappa$. 
One should be able to give upper and lower bounds for $\kappa$, but we have made no attempt
in this direction.

Theorem \ref{main} shows that the volume measure $\mathrm{Vol}$
is completely determined by the metric on $\bm_\infty$. Although this
result seemed plausible, it was not obvious from the construction of the
Brownian sphere in terms of Brownian motion indexed by the Brownian tree.
We note that other models of random geometry such as the Brownian plane
and the Brownian disk have been investigated in recent papers (see in particular
\cite{BM,Plane,CLG,Disks}) and also correspond to specific quantum surfaces,
in the terminology of \cite{MS2,MS3,MS4}
(see \cite[Corollary 1.5]{MS3}). It is not hard to verify that  Theorem \ref{main} 
can be extended to these models, using the known connections between them and the Brownian sphere. 

Let us briefly comment on the proof of Theorem \ref{main}. As is often 
the case in the evaluation of Hausdorff measures, a key ingredient
consists in finding good estimates for the volume of balls. We concentrate
on balls centered at the distinguished point $\bx_*$ of the Brownian sphere,
but the re-rooting invariance property (see Proposition \ref{re-rooting} below) ensures
that similar estimates hold for balls centered at a ``typical point'',
meaning a point chosen uniformly according to the volume measure.
In the construction of the Brownian sphere from the 
Brownian snake driven by a Brownian excursion, $\bx_*$
corresponds to the point with minimal spatial position.
Writing $B(a,r)$ for the closed ball of radius $r$ centered at the point $a$ of $\bm_\infty$,
we are able to show that the $p$-th moment of $\mathrm{Vol}(B(\bx_*,r))$
is bounded above by $C_0^pp!\,r^{4p}$, where $C_0$ is a constant (Proposition 
\ref{moment-normalized}). This bound sharpens a weaker estimate derived 
in \cite[Lemma 6.1]{Acta}. The proof relies on a careful analysis based on a formula 
of \cite{LGW}. Interestingly, estimates for the moments of the volume of balls are discussed in the 
more general setting of the Liouville quantum gravity metric \cite{DDDF,GM} 
in the recent work of Ang, Falconet and Sun \cite{AFS}.

From our estimates on moments of $\mathrm{Vol}(B(\bx_*,r))$ and the re-rooting
invariance property, it is easy to obtain the existence of a constant $K_1$
such that
$$\limsup_{r\downarrow  0} \frac{\mathrm{Vol}(B(a,r))}{h(r)} \leq K_1,$$
for $\mathrm{Vol}$-almost every $a\in\bm_\infty$. Then classical comparison 
results for Hausdorff measures (see Lemma \ref{compa} below) allow us to
find another 
constant $\kappa_1>0$ such that  $m^h(A)\geq \kappa_1\,\mathrm{Vol}(A)$ 
for every Borel subset $A$ of $\bm_\infty$ (Proposition \ref{lower-Haus}). 
In order to get a corresponding upper bound for the Hausdorff measure, we
may again rely on the comparison results, but we need to verify that
the $h$-Hausdorff measure of the set of all points $a\in\bm_\infty$
such that 
$$\limsup_{r\downarrow  0} \frac{\mathrm{Vol}(B(a,r))}{h(r)} < K_2$$
is zero provided that the constant $K_2$ is small enough (Proposition \ref{upper-Haus}).  
This is the most delicate technical part of the paper. We use a ``spine decomposition'' of the
Brownian snake conditioned on its minimal value: This spine decomposition 
provides enough independence between the volumes $\mathrm{Vol}(B(\bx_*,2^{-k}))$
when $k$ varies in $\N$ so that we can bound the probability that these 
volumes are simultaneously small for all $k\in\{k_0,k_0+1,\ldots,n\}$, and we then rely
on the re-rooting invariance property. 

Once we know that the $h$-Hausdorff measure 
is bounded above and below by a positive constant times the volume measure,
we need a kind of zero-one law argument to get that we have indeed 
$m^h=\kappa\,\mathrm{Vol}$ for some constant $\kappa$. Here, the idea is
to consider the canonical projection $\bp$ from $[0,1]$ onto $\bm_\infty$,
and the function $t\mapsto m^h(\bp([0,t]))$. From the bounds on $m^h$,
we know that this function is absolutely continuous, and the point is to prove
that its (almost everywhere defined) derivative is equal to a constant.
This will follow if we can verify that
$$\limsup_{\ve\to 0} \frac{m^h(\bp([s-\ve,s+\ve]))}{2\ve}$$
is equal to a constant $\kappa$, for Lebesgue almost every $s\in(0,1)$
(Lemma \ref{zero-one}). To get this last property, it is convenient to consider the 
free Brownian sphere, whose construction is based on the Brownian snake driven 
by a Brownian excursion distributed according to the It\^o measure. This makes it
possible to use the Bismut decomposition of the Brownian excursion
at a time $U$ chosen according to Lebesgue measure on its duration interval. The technical part of the proof
is to verify that the limsup in the last display (with $s$ replaced by $U$) is
measurable with respect to an appropriate asymptotic $\sigma$-field containing
only events of probability zero or one.

The paper is organized as follows. Section \ref{sec:preli} contains a number of
preliminaries. We recall the construction of the Brownian sphere, and the spine decomposition
of the Brownian snake conditioned on its minimum. We also state the comparison lemma
for Hausdorff measures that plays a central role in our proofs. Section \ref{sec:esti} proves
our estimates on moments of the volume of balls centered at $\bx_*$, from which it
is relatively easy to derive the lower bound  $m^h(A)\geq \kappa_1\,\mathrm{Vol}(A)$.
The proof of the corresponding upper bound is given in Section \ref{sec:upper}. Finally,
the zero-one law argument needed to establish our main result is presented
in Section \ref{proof-main}.

\medskip\noindent
{\it Acknowledgment.} I thank two anonymous referees for several useful remarks on the
first version of this work.

\section{Preliminaries}
\label{sec:preli}

Our main goal in this section is to recall the construction of the Brownian sphere
from the Brownian snake excursion measure. We start with a brief 
discussion of snake trajectories.

\subsection{Snake trajectories}
\label{sna-tra}

By definition,
a finite path $\w$ is a continuous mapping $\w:[0,\zeta]\la\R$, where the
number $\zeta=\zeta_{(\w)}\geq 0$ is called the lifetime of $\w$. We let 
$\W$ denote the space of all finite paths, which is a Polish space when equipped with the
distance
$$d_\W(\w,\w')=|\zeta_{(\w)}-\zeta_{(\w')}|+\sup_{t\geq 0}|\w(t\wedge
\zeta_{(\w)})-\w'(t\wedge\zeta_{(\w')})|.$$
The endpoint or tip of the path $\w$ is denoted by $\wh \w=\w(\zeta_{(\w)})$.
For $x\in\R$, we
set $\W_x=\{\w\in\W:\w(0)=x\}$. The trivial element of $\W_x$ 
with zero lifetime is identified with the point $x$ of $\R$.
%%Occasionally we will use the notation $\underline\w=\min\{\w(t):0\leq t\leq \zeta_{(\w)}\}$.

\begin{definition}
\label{def:snakepaths}
Let $x\in\R$. 
A snake trajectory with initial point $x$ is a continuous mapping $s\mapsto \omega_s$
from $\R_+$ into $\W_x$ 
which satisfies the following two properties:
\begin{enumerate}
\item[\rm(i)] We have $\omega_0=x$ and the number $\sigma(\omega):=\sup\{s\geq 0: \omega_s\not =x\}$,
called the duration of the snake trajectory $\omega$,
is finite (by convention $\sigma(\omega)=0$ if $\omega_s=x$ for every $s\geq 0$). 
\item[\rm(ii)] {\rm (Snake property)} For every $0\leq s\leq s'$, we have
$\omega_s(t)=\omega_{s'}(t)$ for every $t\in[0,\displaystyle{\min_{s\leq r\leq s'}} \zeta_{(\omega_r)}]$.
\end{enumerate} 
\end{definition}

We will write $\S_x$ for the set of all snake trajectories with initial point $x$
and $\S=\bigcup_{x\in\R}\S_x$ for the set of all snake trajectories. If $\omega\in \S$, it is convenient to write $W_s(\omega)=\omega_s$ and $\zeta_s(\omega)=\zeta_{(\omega_s)}$
for every $s\geq 0$ (and we often omit $\omega$ in the notation if there is no ambiguity). The set $\S$ is a Polish space for the distance
$d_{\S}(\omega,\omega')= |\sigma(\omega)-\sigma(\omega')|+ \sup_{s\geq 0} \,d_\W(W_s(\omega),W_{s}(\omega'))$.
A snake trajectory $\omega$ is completely determined 
by the knowledge of the lifetime function $s\mapsto \zeta_s(\omega)$ and of the tip function $s\mapsto \wh W_s(\omega)$: See \cite[Proposition 8]{ALG}.
The range of $\omega$ is defined by
$$\mathcal{R}(\omega):=\{ \wh W_s(\omega):0\leq s\leq \sigma\}$$%=\{W_s(t):0\leq t\leq \zeta_s,0\leq s\leq \sigma\},$$
and we will also use the notation $W_*(\omega):=\min \mathcal{R}(\omega)$. We sometimes write $\omega_*$ instead of $W_*(\omega)$.

Let $\omega\in \S$ be a snake trajectory and $\sigma=\sigma(\omega)$. The lifetime function $s\mapsto \zeta_s(\omega)$ codes a
compact $\R$-tree, which will be denoted 
by $\t_{(\omega)}$ and called the {\it genealogical tree} of the snake trajectory. If $\sim_{(\omega)}$ denotes the equivalence relation
on $[0,\sigma]$ defined by
$$s\sim_{(\omega)} s'\ \hbox{if and only if }\ \zeta_s(\omega)=\zeta_{s'}(\omega)= \min_{s\wedge s'\leq r\leq s\vee s'} \zeta_r(\omega),$$
the $\R$-tree $\t_{(\omega)}$ is the quotient space $[0,\sigma]/\!\sim_{(\omega)}$,
which is equipped with the distance induced by
$$d_{(\omega)}(s,s')= \zeta_s(\omega)+\zeta_{s'}(\omega)-2 \min_{s\wedge s'\leq r\leq s\vee s'} \zeta_r(\omega)$$
(notice that $d_{(\omega)}(s,s')=0$ if and only if $s\sim_{(\omega)} s'$).
We write $p_{(\omega)}:[0,\sigma]\la \t_{(\omega)}$
for the canonical projection, and we root the tree $\t_{(\omega)}$ at $p_{(\omega)}(0)=p_{(\omega)}(\sigma)$. From the snake property, it is immediate that $W_s(\omega)=W_{s'}(\omega)$ 
if $p_{(\omega)}(s)=p_{(\omega}(s')$, so that the mapping $s\mapsto W_s(\omega)$ may be viewed 
as defined on $\t_{(\omega)}$. We sometimes call $\wh W_s(\omega)$ the label of the ``vertex'' $p_{(\omega)}(s)$ of $\t_{(\omega)}$.
%The volume measure on $\t_{(\omega)}$ is defined as the pushforward of
%Lebesgue measure on $[0,\sigma]$ under $p_{(\omega)}$. 

We finally introduce the re-rooting operation on snake trajectories (see \cite[Section 2.2]{ALG}). Let $\omega\in \S_0$ and
$r\in[0,\sigma(\omega)]$. Then $\omega^{[r]}$ is the snake trajectory in $\S_0$ such that
$\sigma(\omega^{[r]})=\sigma(\omega)$ and for every $s\in [0,\sigma(\omega)]$,
\begin{align*}
\zeta_s(\omega^{[r]})&= d_{(\omega)}(r,r\oplus s),\\
\wh W_s(\omega^{[r]})&= \wh W_{r\oplus s}(\omega)-\wh W_r(\omega),
\end{align*}
where we use the notation $r\oplus s=r+s$ if $r+s\leq \sigma(\omega)$, and $r\oplus s=r+s-\sigma(\omega)$ otherwise. 
These prescriptions completely determine $\omega^{[r]}$.
The genealogical tree $\t_{(\omega^{[r]})}$ may be identified to the tree $\t_{(\omega)}$ re-rooted at the vertex $p_{(\omega)}(r)$
\cite[Lemma 2.2]{DLG} (in this identification, the point $p_{(\omega^{[r]})}(s)$ of $\t_{(\omega^{[r]})}$ corresponds to
the point $p_{(\omega)}(r\oplus s)$ of $\t_{(\omega)}$). 

\subsection{The Brownian snake excursion measure}
\label{sec:Excu-mea}

Let $x\in\R$. The Brownian snake excursion 
measure $\N_x$ is the $\sigma$-finite measure on $\S_x$ that satisfies the following two properties: Under $\N_x$,
\begin{enumerate}
\item[(i)] the distribution of the lifetime function $(\zeta_s)_{s\geq 0}$ is the It\^o 
measure of positive excursions of linear Brownian motion, normalized so that the
density of $\sigma$ under $\N_x$ is the function $t\mapsto(2\sqrt{2\pi t^3})^{-1}$,
\item[(ii)] conditionally on $(\zeta_s)_{s\geq 0}$, the tip function $(\wh W_s)_{s\geq 0}$ is
a Gaussian process with mean $x$ and covariance function 
$$K(s,s'):= \min_{s\wedge s'\leq r\leq s\vee s'} \zeta_r.$$
\end{enumerate}
%Informally, the lifetime process $(\zeta_s)_{s\geq 0}$ evolves under $\N_x$ like a Brownian excursion,
%and conditionally on $(\zeta_s)_{s\geq 0}$, each path $W_s$ is a linear Brownian path started from $x$ with lifetime $\zeta_s$, which
%is ``erased'' from its tip when $\zeta_s$ decreases and is ``extended'' when $\zeta_s$ increases. We note that the
%density of $\sigma$ under $\N_0$ is $(2\sqrt{2\pi s^3})^{-1}$. 
The measure $\N_x$ can be interpreted as the excursion measure away from $x$ for the 
Markov process in $\W_x$ called the (one-dimensional) Brownian snake.
We refer to 
\cite{Zurich} for a detailed study of the Brownian snake. For every $t>0$, we can also consider the
conditional probability measure $\N_x^{(t)}:=\N_x(\cdot\midd \sigma =t)$. In fact, $\N_x^{(t)}$
may be defined by the same properties (i) and (ii), just replacing the It\^o measure in (i)
by the law of a positive Brownian excursion with duration $t$. We note that
\begin{equation}
\label{decomp-N}
\N_x=\int_0^\infty \frac{\dd t}{2\sqrt{2\pi t^3}}\,\N^{(t)}_x.
\end{equation}
If $s\in[0,t]$, 
$\N_0^{(t)}$ is invariant under the re-rooting operation $\omega\mapsto \omega^{[s]}$ (see e.g. \cite[Theorem 2.3]{LGW}).

For every $y<x$, we have
\begin{equation}
\label{min-snake}
\N_x(W_*\leq y)=\N_x(y\in\mathcal{R})={\displaystyle \frac{3}{2(x-y)^2}}.
\end{equation}
See e.g. \cite[Section VI.1]{Zurich} for a proof. Additionally, one can prove that $\N_x(\dd \omega)$ a.e., or $\N^{(t)}_x(\dd \omega)$ a.e, there is
a unique $s_*\in[0,\sigma]$ such that $\wh W_{s_*}=W_*$ (see \cite[Proposition 2.5]{LGW}).

The following scaling property is often useful. For $\lambda>0$, for every 
$\omega\in \S_x$, we define $\Theta_\lambda(\omega)\in \S_{x\sqrt{\lambda}}$
by $\Theta_\lambda(\omega)=\omega'$, with
$$\omega'_s(t):= \sqrt{\lambda}\,\omega_{s/\lambda^2}(t/\lambda)\;,\quad
\hbox{for } s\geq 0\hbox{ and }0\leq t\leq \zeta'_s:=\lambda\zeta_{s/\lambda^2}.$$
Then it is a simple exercise to verify that the pushforward of $\N_x$ under $\Theta_\lambda$ is  $\lambda\, \N_{x\sqrt{\lambda}}$. 
Moreover, for every $t>0$, the pushforward of $\N_x^{(t)}$ under $\Theta_\lambda$ is $\N_{x\sqrt{\lambda}}^{(\lambda^2t)}$. 

\subsection{The Brownian sphere}
\label{Brown-sphere}

Let us fix a snake trajectory $\omega\in\S_0$ with 
duration $\sigma=\sigma(\omega)$. For every $s,t\in[0,\sigma]$, we make the 
convention that $[s,t]=[s,\sigma]\cup [0,t]$ if $s>t$ (and of course, if $s\leq t$,
$[s,t]$ is the usual interval). 
We then set
\begin{equation}
\label{formula-D-0}
D^\circ_{(\omega)}(s,t):=\wh W_s(\omega) + \wh W_t(\omega) - 2 \max\Big(\min_{r\in[s,t]} \wh W_r(\omega),\min_{r\in[t,s]}\wh W_r(\omega)\Big).
\end{equation}
and 
\begin{equation}
\label{formula-D}
D_{(\omega)}(s,t)=\inf\Big\{ \sum_{i=1}^p D_{(\omega)}^\circ(t_{i-1},s_i)\Big\},
\end{equation}
where the infimum is over all choices of the integer $p\geq 1$ and of the 
reals $s_0,t_0,s_1,t_1,\ldots,s_p,t_p$ in $[0,\sigma]$ such that $s_0=s$, $t_p=t$
and $p_{(\omega)}(s_i)=p_{(\omega)}(t_i)$ for every $i\in\{0,1,\ldots,p\}$. 
Obviously $D_{(\omega)}\leq D_{(\omega)}^\circ$ (take $p=1$, $t_0=s$ and $s_1=t$).

Clearly, we have $D_{(\omega)}^\circ(s,t)\geq |\wh W_s-\wh W_t|$ for every $s,t\in[0,\sigma]$, and it
follows that we have also $D_{(\omega)}(s,t)\geq |\wh W_s-\wh W_t|$ (recall that 
$p_{(\omega)}(s_i)=p_{(\omega)}(t_i)$ implies $W_{s_i}=W_{t_{i}}$). Additionally, the following simple fact
will be useful in the proof of one of the subsequent technical lemmas. Suppose that
$s_0,t_0,\ldots,s_p,t_p$ are as in formula \eqref{formula-D}. Then, using the
bound $D_{(\omega)}^\circ(s,t)\geq |\wh W_s-\wh W_t|$ and the triangle inequality,
one easily obtains that, for every $j\in\{1,\ldots,p\}$,
\begin{equation}
\label{formula-D-tech}
\sum_{i=1}^p D_{(\omega)}^\circ(t_{i-1},s_i)\geq \wh W_s(\omega) + \wh W_t(\omega) - 2 \max\Big(\min_{r\in[t_{j-1},s_j]} \wh W_r(\omega),\min_{r\in[s_j,t_{j-1}]}\wh W_r(\omega)\Big).
\end{equation}

\rem Let $r\in[0,\sigma]$ be such that $\wh W_r=W_*$. Then, for every $t\in[0,\sigma]$, we have trivially 
$D_{(\omega)}^\circ(r,t)=\wh W_t- W_*=\wh W_t-\wh W_{r}$, and, since we already know that $D^\circ_{(\omega)}(r,t)\geq D_{(\omega)}(r,t)\geq |\wh W_t-\wh W_r|$, we conclude that $D_{(\omega)}(r,t)=\wh W_t- W_*$. 

\smallskip
The mapping $(s,t)\mapsto D_{(\omega)}(s,t)$ defines a pseudo-distance on $[0,\sigma]$, and
we may consider the associated equivalence relation
$$s\approx_{(\omega)} t\hbox{ if and only if }D_{(\omega)}(s,t)=0.$$
Then $D_{(\omega)}$ induces a distance on the quotient space $[0,\sigma]/\!\approx_{(\omega)}$, and we write 
$\bp_{(\omega)}$ for the canonical projection from $[0,\sigma]$ onto $[0,\sigma]/\!\approx_{(\omega)}$. 

The preceding
considerations apply to a fixed snake trajectory $\omega$, but we now randomize 
$\omega$ in order to construct the Brownian sphere. To simplify notation, we usually write $\bp$, $D$ and $D^\circ$
instead of $\bp_{(\omega)}$, $D_{(\omega)}$ and $D^\circ_{(\omega)}$. Recall that
$s_*$ is defined $\N^{(1)}_0$ a.e. as the unique element of $[0,1]$ such that $\wh W_{s_*}=W_*$.

\begin{definition}
\label{def:sphere}
The standard Brownian sphere is defined under the probability measure $\N^{(1)}_0(\dd \omega)$ as the random measure metric space $\bm_\infty:=[0,1]/\!\approx_{(\omega)}$ equipped 
with the distance induced by $D$ (for which we keep the same notation $D$) and with the volume measure $\mathrm{Vol}$ which is the pushforward of Lebesgue measure on $[0,1]$ under the
canonical projection $\bp$.%, and with the distinguished point $\mathbf{x}_*:=\bp(s_*)$. 
\end{definition}

We could also have introduced the free Brownian sphere, which is defined in the same way replacing 
$\N^{(1)}_0$ by $\N_0$ (and the interval $[0,1]$ by $[0,\sigma]$). As a side remark, we note that
the property $s\sim_{(\omega)} t$ obviously implies $s\approx_{(\omega)} t$, and thus one may as well define $\bm_\infty$
as a quotient space of the tree $\t_{(\omega)}$ under $\N^{(1)}_0$ (this is the point of view of \cite{Uniqueness} in particular).
%One sometimes defines the  Brownian sphere as a $2$-pointed measure metric space,
%letting $\bx_0:=\bp(0)$ be the second distinguished point. 

In the present work, we will view the (standard) Brownian sphere as a {\it pointed} measure metric space, with
the distinguished point $\mathbf{x}_*:=\bp(s_*)$.
So the Brownian sphere is  a random 
variable $(\bm_\infty, D,\mathrm{Vol},\mathbf{x}_*)$ with values in the space $\M^\bullet$ of all isometry classes
of pointed compact measure metric spaces, which is equipped with the
Gromov-Hausdorff-Prokhorov topology. We refer to \cite[Section 2.1]{Disks} for a brief presentation of
the Gromov-Hausdorff-Prokhorov topology on $\M^\bullet$.

Since $\wh W_{s_*}=W_*$, the remark after \eqref{formula-D-tech} shows that we have
$\N^{(1)}_0$ a.s., for every $t\in [0,1]$,
\begin{equation}
\label{key-dist}
D(s_*,t)=\wh W_t -W_*,
\end{equation}
and consequently $D(\bx_*,a)=\wh W_t -W_*$, for any $a=\bp(t)\in \bm_\infty$.
This property will be crucial
for our applications. 

We now observe that the distinguished point $\mathbf{x}_*$
is not a special point of the Brownian sphere, in the sense that it could be replaced by another point distributed according 
to the volume measure $\mathrm{Vol}$, without changing the distribution of the $4$-tuple $(\bm_\infty, D,\mathrm{Vol},\mathbf{x}_*)$. We state a slightly more precise result.

\begin{proposition}
\label{re-rooting}
Let $t\in[0,1]$. Then $(\bm_\infty, D,\mathrm{Vol},\bp(t))$ has the same distribution as $(\bm_\infty, D,\mathrm{Vol},\mathbf{x}_*)$.
\end{proposition}

\proof The fact that the distribution of $(\bm_\infty, D,\mathrm{Vol},\bp(t))$ does not depend on $t$ is an easy consequence
of the invariance of $\N^{(1)}_0$ under the re-rooting operation (replacing $\omega$ by $\omega^{[t]}$ gives a measure metric space
$\bm_\infty(\omega^{[t]})$ which is the same as $\bm_\infty(\omega)$ modulo a measure-preserving isometry, but 
the point $\bp_{(\omega^{[t]})}(0)$ of $\bm_\infty(\omega^{[t]})$ corresponds to the point $\bp_{(\omega)}(t)$ of $\bm_\infty(\omega)$). The fact that this
distribution is the same as the
distribution of $(\bm_\infty, D,\mathrm{Vol},\mathbf{x}_*)$ follows from \cite[Theorem 8.1]{Acta}.
Note that \cite[Theorem 8.1]{Acta} deals with metric spaces and not with measure metric spaces, and thus one needs
a slight extension of this result, which is however
derived by the very same arguments as in \cite{Acta}, using the convergence of 
rescaled quadrangulations to the Brownian sphere in the Gromov-Hausdorff-Prokhorov sense,
as stated in \cite[Theorem 7]{Disks}. \endproof

\subsection{Decomposing the Brownian snake at its minimum}
\label{deco-mini}

As an important ingredient of our proofs, we will use the conditional distribution of the
Brownian snake under $\N_0$ given its minimum $W_*$. This conditional 
distribution is described in \cite{Bessel} via a ``spine decomposition'' involving
a nine-dimensional Bessel process and two Poisson point measures on the 
space of snake trajectories. We do not give details of this spine decomposition (see 
\cite{Bessel}) but we present the consequence that will
be relevant to the present work.

The results of \cite{Bessel} allow us to make sense of the conditional distribution
$\N_0(\dd\omega\midd W_*=-x)$ for every $x>0$ (in such a way that it depends
continuously on $x$). Via the obvious translation, we can also make
sense of $\N_x(\dd\omega\midd W_*=0)$ for every $x>0$.
% By scaling arguments, the unique time $s_*\in[0,\sigma]$
%minimizing $\wh W_s$ is also well defined under $\N_x(\dd\omega\midd W_*=0)$.

\begin{proposition}
\label{spine-dec}
Let $x>0$. Almost surely under $\N_x(\dd\omega\midd W_*=0)$, 
we can define a random finite path $(U_t)_{0\leq t\leq L_x}$
and a point measure 
$$\mathcal{N}(\dd t\,\dd\omega')=\sum_{i\in I} \delta_{(t_i,\omega_i)}(\dd t\,\dd\omega'),$$
on $\R_+\times\S$, such that, for every nonnegative measurable function $\varphi$ on $\R$,
$$\int_0^{\sigma(\omega)} \varphi(\wh W_s(\omega))\,\dd s =\sum_{i\in I} \int_0^{\sigma(\omega_i)} \varphi(\wh W_s(\omega_i))\,\dd s,$$
and furthermore the following properties hold under $\N_x(\dd\omega\midd W_*=0)$:
\begin{description}
\item{$\bullet$} $U=(U_t)_{0\leq t\leq L_x}$ 
is distributed as a nine-dimensional Bessel process started from $0$ up to its last passage time at level $x$;
\item{$\bullet$} conditionally on $U$, $\mathcal{N}(\dd t\,\dd\omega')$  is a Poisson point measure on $\R_+\times \S$ with intensity
$$4\,\mathbf{1}_{\{t\leq L_x\}}\,\dd t\,\N_{U_t}(\dd\omega' \cap\{W_*(\omega')>0\}).$$
\end{description}
%and are such that
%$$\mathcal{R}(\omega)= \{U_t: 0\leq t\leq L_x\} \cup \Big(\bigcup_{i\in I} \mathcal{R}(\omega_i)\Big).$$
\end{proposition}

\rem The finite path $U$ is constructed as the time-reversal of the path $W_{s_*}$ 
and the point measure $\n$ accounts for the subtrees branching off the ancestral line of $p_{(\omega)}(s_*)$
in the tree $\t_{(\omega)}$. 

\smallskip
The occurence of the nine-dimensional Bessel process in the preceding spine decomposition is
closely related to an absolute continuity relation between the laws of Bessel processes,
which we now recall in the special case that will be of interest to us. Let $x>0$, and suppose that,
under the probability measure $\P_x$, we are given two processes $(B_s)_{s\geq 0}$ and
$(R_s)_{s\geq 0}$ which are respectively a linear Brownian motion started at $x$ and a 
nine-dimensional Bessel process started at $x$. Then, for every $t>0$ and for every
nonnegative measurable function $F$ on the space $C([0,t],\R)$
of all continuous functions from $[0,t]$ into $\R$, we have
\begin{equation}
\label{Bessel-ac}
\E_x\Big[ \mathbf{1}_{\{B_s>0,\forall s\in[0,t]\}}\,\exp\Big(-6\int_0^t \frac{\dd s}{B_s^2}\Big)\, F\big((B_s)_{0\leq s\leq t}\big)\Big]
= x^4\,\E_x\Big[ R_t^{-4}\,F\big((R_s)_{0\leq s\leq t}\big)\Big].
\end{equation}
See e.g. Exercise XI.1.22 in \cite{RY}. 

\subsection{Hausdorff measures}

Let $E$ be a compact metric space, let $\delta>0$ and let $h$ be a continuous nondecreasing function from $[0,\delta]$ into $\R_+$
such that $h(0)=0$ and $h(r)>0$ for every $r\in(0,\delta]$. If $A$ is a subset of $E$, we denote the diameter of $A$
by $\mathrm{diam}(A)$. Then the Hausdorff measure $m^h(A)\in[0,\infty]$ is defined by
\begin{equation}
\label{def-Haus}
m^h(A)=\lim_{\ve \downarrow 0}  \Bigg(\inf_{(U_i)_{i\in I}\in \mathrm{Cov}_\ve(A)} \sum_{i\in I} h\big(\mathrm{diam}(U_i)\big)\Bigg),
\end{equation}
where $\mathrm{Cov}_\ve(A)$
is the collection of all countable coverings of $A$ by subsets of $E$ with diameter smaller than $\ve$. Note that the quantity inside the
big parentheses in \eqref{def-Haus} is a nonincreasing function of $\ve$, so that the limit always exists in $[0,\infty]$. 

As a key ingredient of the proof of our main results, we will use certain comparison 
results for Hausdorff measures, which we now recall for the reader's convenience.
The following statement can be found in \cite[Lemma 2.1]{LGD}, and is a mild 
generalization of results proved by Rogers and Taylor \cite{RT} for
subsets of Euclidean space. If $x\in E$ and $r>0$, we denote the closed ball of 
radius $r$  centered at $x$ by $\mathcal{B}(x,r)$.

\begin{lemma}
\label{compa}
Assume that the function $h$ satisfies $h(2r)\leq ch(r)$ for every $r\in(0,\delta/2]$, for some constant $c>1$, and
let $\mu$ be a finite Borel measure on $E$. There exist two
positive constants $M_1$ and $M_2$, which only depend on $c$, such that the following holds for every Borel subset $A$ of 
$E$ and every $b>0$.
\begin{itemize}
\item[\rm(i)] If, for every $x\in A$,
$$\limsup_{r\downarrow 0} \frac{\mu(\mathcal{B}(x,r))}{h(r)} \leq b,$$
then $m^h(A)\geq M_1b^{-1}\,\mu(A)$.
\item[\rm(ii)] If, for every $x\in A$,
$$\limsup_{r\downarrow 0} \frac{\mu(\mathcal{B}(x,r))}{h(r)} \geq b,$$
then $m^h(A)\leq M_2b^{-1}\,\mu(A)$.
\end{itemize}
\end{lemma}

We conclude this section with a technical point concerning the measurability of
functions defined as the Hausdorff measure of certain particular subsets of
the Brownian sphere $\bm_\infty$. In Section \ref{proof-main} below, especially in the proof of Lemma \ref{zero-one},
we will consider quantities of the type $m^h(\bp([u,v]))$, where $0\leq u\leq v$,
and it is not immediately obvious from the definition that these quantities are
random variables. This measurability question can be settled as follows. 
We consider a countable dense subset $\mathcal{D}$ of $\bm_\infty$, for instance
the set of all $\bp(r)$ for rational values of $r$, and the collection $\mathcal{C}$ of 
all closed balls of rational radius centered at a point of $\mathcal{D}$.
A compactness argument shows that, in order to evaluate the infima  
appearing in formula \eqref{def-Haus} for $m^h(\bp([u,v]))$, we may restrict our attention to {\it finite}
coverings of $\bp([u,v])$ by subsets of $\bm_\infty$ which are {\it finite}
unions of balls in $\mathcal{C}$. In this way, we see that \eqref{def-Haus} only involves 
infima of countably many random variables, giving the desired measurability property.

\section{Estimates for the volume of balls}
\label{sec:esti}

Our main objective in this section is to prove that, for any Borel subset $A$ of $\bm_\infty$, the Hausdorff measure 
$m^h(A)$ with respect to the gauge function $h$ introduced in Theorem \ref{main} 
is bounded below by a constant times $\mathrm{Vol}(A)$ (see Proposition \ref{lower-Haus} below). 
To this end, we will use part (i) of Lemma \ref{compa}, and we thus need to get
precise bounds on the volume of balls in the Brownian 
sphere $\bm_\infty$, which is defined under the probability measure
$\N_0^{(1)}=\N_0(\cdot\midd \sigma =1)$ introduced in Section \ref{sec:Excu-mea}. Thanks to Proposition \ref{re-rooting}, it will be enough to
consider balls at the distinguished point $\bx_*$. For $a\in\bm_\infty$ and $r>0$, let $B(a,r)$ denote the closed 
ball of radius $r$ centered at $a$ in $\bm_\infty$. Then property \eqref{key-dist} and the definition
of the volume measure give
\begin{equation}
\label{vol-ball-x*}
\mathrm{Vol}(B(\bx_*,r))=\int_0^1 \mathbf{1}_{\{\wh W_s-W_*\leq r\}}\,\dd s.
\end{equation}
Much of this section is thus devoted to the proof of the following proposition.

\begin{proposition}
\label{moment-normalized}
There exists a constant $C_0$ such that, for every integer $p\geq 1$ and every $\ve>0$,
$$\N^{(1)}_0\Big( \Big(\int_0^1 \dd s\,\mathbf{1}_{\{\wh W_s - W_* \leq\ve\}}\Big)^p\Big) \leq C_0^p\, p!\,\ve^{4p}.$$
\end{proposition}

The proof of Proposition \ref{moment-normalized} depends on several intermediate results.
Although we are primarily interested in the
probability measure $\N_0^{(1)}$, it will be useful to derive first certain preliminary estimates
under the ($\sigma$-finite) measures $\N_x$.
In particular, Proposition \ref{key-estim} below provides key bounds under $\N_x$, which will also play an 
important role in the next section.  Before stating and proving Proposition \ref{key-estim}, we
need to recall a formula of \cite{LGW} giving the expectation under $\N_0$
of certain functionals of the Brownian snake. We first introduce the relevant notation.
We follow closely
\cite[Section 2.2]{LGW}. 

Let 
$$\mathcal{U}:=\bigcup_{n=0}^\infty \{1,2\}^n$$
where by convention $\{1,2\}^0=\{\varnothing\}$. An element of $\mathcal{U}$
is a finite sequence $u=u^1\ldots u^n$ of elements of $\{1,2\}$, and 
we write $|u|=n$ for the length of this sequence (if $u=\varnothing$, $|u|=0$). The mapping
$\pi:\mathcal{U}\backslash\{\varnothing\}\la\mathcal{U}$ is defined by $\pi(u^1\ldots u^n)=u^1\ldots u^{n-1}$.
For every integer $k\geq 1$, $\pi^k$ stands for the $k$-th iterate of $\pi$ (note that $\pi^k$
is defined on $\{u\in\mathcal{U}:|u|\geq k\}$). If $u=u^1\ldots u^n\in\mathcal{U}$, we will use the
obvious notation $u1=u^1\ldots u^n1$ and $u2=u^1\ldots u^n2$ for the two elements $v$ of $\mathcal{U}$
such that $\pi(v)=u$. 

By definition, a binary (plane) tree is a finite subset $\tau$ of $\mathcal{U}$ which satisfies the following properties.
\begin{enumerate}
\item[(i)] $\varnothing\in \tau$.
\item[(ii)] If $u\in \tau$ and $u\not=\varnothing$,	 then $\pi(u)\in \tau$.
\item[(iii)] For every $u\in\tau$, either $u1\in\tau$ and $u2\in \tau$, or 
$u1\notin \tau$ and $u2\notin \tau$. In the latter case we say that $u$ is a leaf of $\tau$.
\end{enumerate}

We denote the set of all binary trees by $\mathbf{T}^b$. Then a marked (binary) tree is
a pair $(\tau,(\ell_u)_{u\in\tau})$, where $\tau\in\mathbf{T}^b$, and
$\ell_u\in[0,\infty)$, for every $u\in\tau$. It will be convenient to view a marked tree as the compact $\R$-tree obtained
by gluing line segments of length $\ell_u$, for every $u\in\tau$, 
according to the genealogical structure prescribed by $\tau$. To make this more precise,
let $\theta=(\tau,(\ell_u)_{u\in\tau})$ be a marked tree, and
introduce the vector space $\R^\tau$ of all mappings from $\tau$ into $\R$, which is equipped
with the usual topology. Write 
$(\ve_u)_{u\in\tau}$ for the canonical basis of $\R^\tau$. Then consider the mapping
$$p_\theta:\bigcup_{u\in\tau} \{u\}\times [0,\ell_u] \la \R^\tau,$$
defined by
$$p_\theta(u,\ell)=\sum_{k=1}^{|u|} \ell_{\pi^k(u)}\ve_{\pi^k(u)} + \ell\,\ve_u.$$
The $\R$-tree $\tilde\theta$ associated with $\theta$ is the range of $p_\theta$
(which is a connected union of line segments in $\R^\tau$)
equipped with the distance $d_\theta$ such that $d_\theta(a,b)$ 
is the length of the unique (up to reparameterization) continuous injective path
from $a$ to $b$ in $\tilde\theta$. We write ${\mathrm{Leb}}_\theta$ for Lebesgue measure
on $\tilde \theta$: ${\mathrm{Leb}}_\theta$ is obtained as the sum over all
$u\in\tau$ of the pushforward 
of Lebesgue measure on $[0,\ell_u]$ under the mapping $\ell\mapsto p_\theta(u,\ell)$. By definition, leaves of $\tilde\theta$ are the points of
the form $p_\theta(u,\ell_u)$ where $u$ is a leaf of $\tau$. 
On the other hand,
points of the form $p_\theta(u,\ell_u)$ where $u\in\tau$ is not a leaf are called 
nodes of $\tilde\theta$. 
We write $L(\theta)$ for the set
of all leaves of $\tilde \theta$, and $N(\theta)$ for the set of all nodes.

We next introduce Brownian motion indexed by $\tilde \theta$. We fix $x\in\R$,
and, under a probability measure denoted by $Q^\theta_x$, we consider a collection 
$(\xi^u)_{u\in\tau}$ of independent linear Brownian motions, which all start from $0$
except for $\xi^\varnothing$, which starts from $x$. We then define a 
continuous random process $(V_a)_{a\in\tilde\theta}$
by setting
$$V_{p_\theta(u,\ell)}=\sum_{k=1}^{|u|} \xi^{\pi^k(u)}(\ell_{\pi^k(u)}) + \xi^u(\ell),$$
for every $u\in\tau$ and $\ell\in[0,\ell_u]$. 

It will also be convenient to
assume that, under the probability measure $Q^\theta_x$, in addition to the
collection $(\xi^u)_{u\in\tau}$, we are given an independent Poisson point measure on
$\tilde \theta\times \S_0$  with intensity
$4\,{\mathrm{Leb}}_\theta(\dd a)\,\N_0(\dd\omega)$, and we denote this point measure by $\sum_{j\in J}\delta_{(a_j,\omega_j)}$. 

Fix an integer $p\geq 1$. Let $\mathbf{T}^b_p$ stand for the set of all 
binary trees with $p$ leaves, and let $\T_p$ be the corresponding set of 
marked trees. We define Lebesgue measure $\Lambda_p$ on $\T_p$
by the formula
$$\int_{\T_p} \Lambda_p(\dd\theta)\,F(\theta) = \sum_{\tau\in\mathbf{T}^b_p} \int \prod_{u\in\tau} \dd \ell_u\,F\Big(\tau,(\ell_u)_{u\in\tau}\Big).$$

Finally, we write $\mathcal{K}$ for the set of all compact subsets of $\R$, which is equipped
with the usual Hausdorff metric. Recall our notation $\mathcal{R}(\omega)=\{\wh W_s(\omega):0\leq s\leq \sigma\}$ for the range of 
$\omega$. The next proposition is a special case of \cite[Theorem 2.2]{LGW} (see also the remark after this theorem).

\begin{proposition}
\label{formula-p-moment}
 Let $F$ be a nonnegative measurable function on $\R^p\times \mathcal{K}$, which is symmetric
 with respect to the coordinates of $\R^p$. Then
\begin{align*}
&\N_0\Bigg(\int_{[0,\sigma]^p} \dd s_1\ldots, \dd s_p\,F\big((\wh W_{s_1},\ldots,\wh W_{s_p}),\mathcal{R}\big)\Bigg)\\
 &\qquad= 2^{p-1}p!\int \Lambda_p(\dd\theta)\, Q^\theta_x\Bigg[ F\Bigg( (V_a)_{a\in L(\theta)}, \{V_a:a\in\tilde\theta\}\cup\Big(\bigcup_{j\in J}(V_{a_j}+\mathcal{R}(\omega_j))
 \Big)\Bigg)\Bigg].
 \end{align*}
%where $\mathrm{cl}(A)$ denotes the closure of a subset $A$ of $\R$.
 \end{proposition}
 
 Proposition \ref{formula-p-moment} is a crucial ingredient in the proof of the following key estimate.
 
 \begin{proposition}
 \label{key-estim}
 There exist two positive constants $C_1$ and $C_2$ such that,
 for every $\ve,x>0$,
 $$C_1\,x^2\,\ve^{5}\,(x\vee \ve)^{-5}\leq \N_x\Bigg( \Big(\int_0^\sigma \dd s\,\mathbf{1}_{\{\wh W_s\leq \ve\}}\Big)\,\mathbf{1}_{\{W_*>0\}}\Bigg)
 \leq C_2\,x^{2}\,\ve^{5}\,(x\vee \ve)^{-5},$$
 and for every integer $p\geq 2$,
 $$C_1^p\,p!\,x^4\,\ve^{1+4p}\,(x\vee \ve)^{-7}\leq \N_x\Bigg( \Big(\int_0^\sigma \dd s\,\mathbf{1}_{\{\wh W_s\leq \ve\}}\Big)^p\,\mathbf{1}_{\{W_*>0\}}\Bigg)
 \leq C_2^p\,p!\,x^{4}\,\ve^{1+4p}\,(x\vee \ve)^{-7}.$$
 \end{proposition}
 
 \proof Let $x>0$ and $\ve>0$, and let $p\geq 1$ be an integer.
We apply Proposition \ref{formula-p-moment} with 
$$F\big((\wh W_{s_1},\ldots,\wh W_{s_p}),\mathcal{R}\big)= \mathbf{1}_{\{0<\wh W_{s_1}\leq\ve,\ldots,0<\wh W_{s_p}\leq\ve\}}\, \mathbf{1}_{\{\mathcal{R}\subset(0,\infty)\}}.$$
We get
 \begin{align}
 \label{key-tech1}
 &\N_x\Bigg( \Big(\int_0^\sigma \dd s\,\mathbf{1}_{\{\wh W_s\leq \ve\}}\Big)^p\,\mathbf{1}_{\{W_*>0\}}\Bigg)\nonumber\\
 &\quad= 2^{p-1}p!\int \Lambda_p(\dd\theta)\, Q^\theta_x\Bigg[ \Bigg( \prod_{a\in L(\theta)} \mathbf{1}_{\{0<V_a\leq \ve\}}\Bigg)\,
 \mathbf{1}_{\{V_a>0,\forall a\in\tilde\theta\}}\,\exp\Big(-4 \int {\mathrm{Leb}}_\theta(\dd a)\,\N_{V_a}(W_*\leq 0)\Big)\Bigg].
 \end{align}
 The appearance of the exponential term in the right-hand side of \eqref{key-tech1} comes from the 
 fact that
 $$Q^\theta_x\Bigg[ \Big(\bigcup_{j\in J}(V_{a_j}+\mathcal{R}(\omega_j))
 \Big) \subset (0,\infty) \,\Bigg| \, (V_a)_{a\in\tilde\theta}\Bigg]= \exp\Big(-4 \int {\mathrm{Leb}}_\theta(\dd a)\,\N_{V_a}(W_*\leq 0)\Big)\Bigg],$$
 since $\sum_{j\in J}\delta_{(a_j,\omega_j)}$ is Poisson  with intensity
$4\,{\mathrm{Leb}}_\theta(\dd a)\,\N_0(\dd\omega)$ and is independent of $(V_a)_{a\in\tilde\theta}$ under $Q^\theta_x$. 

 Let us fix $\theta=(\tau,(\ell_u)_{u\in\tau})\in \T_p$. Recalling formula \eqref{min-snake},
 the expectation under $Q^\theta_x$
 in \eqref{key-tech1} is also equal to
 \begin{equation}
 \label{key-tech2}
 I^\ve_x(\theta):= Q^\theta_x\Bigg[ \Bigg( \prod_{a\in L(\theta)} \mathbf{1}_{\{0<V_a\leq \ve\}}\Bigg)\,
 \mathbf{1}_{\{V_a>0,\forall a\in\tilde\theta\}}\,\exp\Big(-6 \int {\mathrm{Leb}}_\theta(\dd a)\,(V_a)^{-2}\Big)\Bigg].
 \end{equation}
 At this stage, we use formula \eqref{Bessel-ac}
 to get rid of the term $\exp(-6 \int {\mathrm{Leb}}_\theta(\dd a)\,(V_a)^{-2})$ in $I^\ve_x(\theta)$. To this end, we introduce a process
$(\ov{V}_a)_{a\in\tilde \theta}$ which is (informally) obtained by running independent nine-dimensional Bessel processes 
along the branches of $\tilde \theta$. To be specific, under the same probability measure $Q^\theta_x$ (for every
$\theta\in \T_p$), we define processes $(\ov\xi^u(t))_{0\leq t\leq \ell_u}$, for every $u\in\tau$, inductively, by first requiring that
$\ov\xi^\varnothing$ is the (unique) solution of the stochastic integral equation
$$\ov\xi^\varnothing(t)=x + \xi^\varnothing(t) + \int_0^t \frac{4}{\ov\xi^\varnothing(s)}\,\dd s\,,\quad 0\leq t\leq \ell_\varnothing$$
and then inductively, for every $u\in\tau\backslash\{\varnothing\}$, $\ov\xi^u$ solves
$$\ov\xi^u(t)=\ov\xi^{\pi(u)}(\ell_{\pi(u)}) + \xi^u(t) + \int_0^t \frac{4}{\ov\xi^u(s)}\,\dd s\,,\quad 0\leq t\leq \ell_u.$$
We then define $(\ov V_a)_{a\in\tilde \theta}$ by the relation $\ov V_{p_\theta(u,\ell)}=\ov\xi^u_\ell$ for every $u\in\tau$ and $\ell\in[0,\ell_u]$. 
Note that we have now $\ov V_a>0$ for every $a\in\tilde \theta$, $Q^\theta_x$ a.s.

We claim that the quantity $I^\ve_x(\theta)$ in \eqref{key-tech2} is also equal to $x^4\,J^\ve_x(\theta)$, where
\begin{equation}
\label{key-tech3}
J^\ve_x(\theta):=Q^\theta_x\Bigg[ \Bigg( \prod_{a\in L(\theta)} \mathbf{1}_{\{\ov V_a\leq \ve\}}\,(\ov V_a)^{-4}\Bigg)\,
\Bigg(\prod_{a\in N(\theta)} (\ov V_a)^4\Bigg)\Bigg]
\end{equation}
and we recall that $N(\theta)$ is the set of all nodes of $\tilde\theta$. 
%This follows from \eqref{Bessel-ac}
%as in the proof of \cite[Theorem 5.1]{LGW}. 
Let us verify our claim by induction on $p$. As in Section \ref{deco-mini}, we consider a
linear Brownian motion $(B_t)_{t\geq 0}$ and a nine-dimensional Bessel process $(R_t)_{t\geq 0}$ that both start at $y$ under the probability measure $\P_y$, for any $y>0$. 
When $p=1$, if $a_1$ denotes the unique leaf of $\tilde \theta$ and $T=\ell_\varnothing$, we have
\begin{align*}
I^\ve_x(\theta)&=Q^\theta_x\Big[\mathbf{1}_{\{V_{a_1}\leq\ve\}}\,
 \mathbf{1}_{\{V_a>0,\forall a\in\tilde\theta\}}\,\exp\Big(-6 \int {\mathrm{Leb}}_\theta(\dd a)\,(V_a)^{-2}\Big)\Big]\\
 &=\E_x\Big[ \mathbf{1}_{\{B_T\leq\ve\}}\,\mathbf{1}_{\{B_t>0,\forall t\in[0,T]\}}\,\exp\Big(-6\int_0^T \dd t\,(B_t) ^{-2}\Big)\Big]\\
&=x^4 \,\E_x\Big[(R_T)^{-4}\,\mathbf{1}_{\{R_T\leq\ve\}}\Big]\\
 &=x^4\,Q^\theta_x\Big[ (\ov V_{a_1})^{-4}\,\mathbf{1}_{\{\ov V_{a_1}\leq \ve\}}\Big]
 \end{align*}
where we used \eqref{Bessel-ac} in the third equality. This gives the case $p=1$ of our claim. Then let $p\geq 2$
and assume that our claim has been proved up to order $p-1$. Let $\theta_1$ and $\theta_2$
be the two marked subtrees that are obtained by ``breaking'' $\theta$ at its first node (more precisely,
$\tilde\theta_1$ and $\tilde \theta_2$ are isometric to the closures of the two connected components 
of $\tilde\theta\,\backslash\, p_\theta(\{\varnothing\}\times [0,\ell_\varnothing])$) and write again $T=\ell_\varnothing$
to simplify notation. By construction, $I^\ve_x(\theta)$
is equal to
\begin{align*}
&\E_x\Bigg[\mathbf{1}_{\{B_t>0,\forall t\in[0,T]\}}\exp\Big(\!\!-6\!\!\int_0^{T}\frac{\dd t}{B_t^2}\Big)\,
Q^{\theta_1}_{B_T}\Big[ \Big( \prod_{a\in L(\theta_1)} \mathbf{1}_{\{V_a\leq \ve\}}\Big)
 \mathbf{1}_{\{V_a>0,\forall a\in\tilde\theta_1\}}\exp\Big(\!\!-6\!\! \int \!{\mathrm{Leb}}_{\theta_1}(\dd a)\,(V_a)^{-2}\Big)\Big]\\
 &\qquad \qquad\qquad\qquad\qquad\qquad\qquad\times Q^{\theta_2}_{B_T}\Big[ \Big( \prod_{a\in L(\theta_2)} \mathbf{1}_{\{V_a\leq \ve\}}\Big)
 \mathbf{1}_{\{V_a>0,\forall a\in\tilde\theta_2\}}\exp\Big(\!\!-6 \!\!\int \!{\mathrm{Leb}}_{\theta_2}(\dd a)\,(V_a)^{-2}\Big)\Big]\Bigg]\\
 &=\E_x\Bigg[\mathbf{1}_{\{B_t>0,\forall t\in[0,T]\}}\exp\Big(\!\!-6\!\!\int_0^{T}\frac{\dd t}{B_t^2}\Big)\,
 (B_{T})^4\,J^\ve_{B_T}(\theta_1) \times (B_T)^4 J^\ve_{B_T}(\theta_2)\Bigg]\\ &=x^4\,\E_x\Big[(R_T)^4 J^\ve_{R_T}(\theta_1)\,J^\ve_{R_T}(\theta_2)\Big],
\end{align*}
using the induction hypothesis in the first equality and again \eqref{Bessel-ac} in the second one. Our claim follows since one immediately verifies
that $J^\ve_x(\theta)= \E_x[(R_T)^4 J^\ve_{R_T}(\theta_1)\,J^\ve_{R_T}(\theta_2)]$.

For every binary tree $\tau\in \mathbf{T}^b_p$, let $\Lambda^{(\tau)}$ be the measure on $\T_p$ defined by
$$\int_{\T_p} \Lambda^{(\tau)}(\dd\theta)\,F(\theta) =\int \prod_{u\in\tau} \dd \ell_u\,F\Big(\tau,(\ell_u)_{u\in\tau}\Big).$$
Note that $\Lambda_p=\sum_{\tau\in \mathbf{T}^b_p}\Lambda^{(\tau)}$, and that $\#\mathbf{T}^b_p$ is the Catalan number of 
order $p-1$ so that $\#\mathbf{T}^b_p\sim \frac{4^{p-1}}{\sqrt{\pi} p^{3/2}}$ as $p\to\infty$. Using these observations 
together with formulas \eqref{key-tech1},\eqref{key-tech2},\eqref{key-tech3}, we see that the statement of
Proposition \ref{key-estim} follows from the next lemma. We write $\tau_0=\{\varnothing\}$ for the unique element of $\mathbf{T}^b_1$.

\begin{lemma}
\label{key-lem}
There exists a positive constant $c_1$ such that, for every $\ve,x>0$,
\begin{equation}
\label{key-lem101}
c_1\,\ve^{5}\,x^{-2}(x\vee \ve)^{-5}\leq \int \Lambda^{(\tau_0)}(\dd \theta)\,J^\ve_x(\theta)\leq \ve^{5}\,x^{-2}\,(x\vee \ve)^{-5},
\end{equation}
and, for every integer $p\geq 2$ and every
$\tau\in\mathbf{T}^b_p$,
\begin{equation}
\label{key-lem102}
c_1^p\,\ve^{1+4p}(x\vee \ve)^{-7}\leq \int \Lambda^{(\tau)}(\dd \theta)\,J^\ve_x(\theta)\leq \ve^{1+4p}(x\vee \ve)^{-7}.
\end{equation}
\end{lemma}

 \proof
We first observe that
$$\int \Lambda^{(\tau_0)}(\dd \theta)\,J^\ve_x(\theta)=\int_0^\infty \dd t\,\E_x[\mathbf{1}_{\{R_t\leq \ve\}}R_t^{-4}].$$
To evaluate this quantity, 
 it is convenient to write $(Z_t)_{t\geq 0}$ for a nine-dimensional Brownian motion that starts from 
 $z$ under the probability measure $\P^{(9)}_z$, for every $z\in \R^9$. We write $G(z,z')=c_0|z-z'|^{-7}$
 for the Green function of $Z$, where $c_0=\frac{1}{2}\pi^{-9/2}\Gamma(7/2)$. Also, for every $x>0$, let $z_x$
 denote a point of $\R^9$ such that $|z_x|=x$, and let $\Sigma(\dd y)$ stand for the volume measure
 on the unit sphere $\mathbb{S}^8$ of $\R^9$. Recall that $\Sigma(\mathbb{S}^8)=2\pi^{9/2}/\Gamma(9/2)$. Then, we have
 \begin{align*}
 \int_0^\infty \dd t\,\E_x[\mathbf{1}_{\{R_t\leq \ve\}}R_t^{-4}]&=\E^{(9)}_{z_x}\Big[\int_0^\infty \dd t\,\mathbf{1}_{\{|Z_t|\leq \ve\}}\,|Z_t|^{-4}\Big]\\
 &=\int_{\{|z|\leq \ve\}} \dd z\,|z|^{-4}\,G(z_x,z)\\
 &=c_0\int_0^\ve \dd \rho\,\rho^4\,\int_{\mathbb{S}^8} \Sigma(\dd y)\,|\rho y-z_x|^{-7}\\
 &=\frac{2}{7}\int_0^\ve \dd \rho\,\rho^4\,(x\vee\rho)^{-7},
 \end{align*}
 using the formula $\int_{\mathbb{S}^8} \Sigma(\dd y)\,|\rho y-z_x|^{-7}=\Sigma(\mathbb{S}^8)\,(x\vee \rho)^{-7}$, which easily
 follows from the fact that the function $z\mapsto |z-z_x|^{-7}$ is harmonic on $\R^9\backslash\{z_x\}$. If $\ve\leq x$,
 we use 
 $$\int_0^\ve \dd\rho\, \rho^4x^{-7}=x^{-7}\,\frac{\ve^5}{5}$$
 and, if $x<\ve$,
 $$\int_0^x \dd \rho\, \rho^4x^{-7} + \int_x^{\ve} \dd\rho \,\rho^{-3}= \frac{x^{-2}}{5} + \frac{1}{2}(x^{-2}-\ve^{-2}).$$
From these elementary calculations, 
we can pick an explicit constant $c_1>0$ such that
 \begin{equation}
\label{key-lem103}
100\,c_1\,\ve^{5}\,x^{-2}(x\vee \ve)^{-5}\leq \int \Lambda^{(\tau_0)}(\dd \theta)\,J^\ve_x(\theta)\leq \ve^{5}\,x^{-2}\,(x\vee \ve)^{-5},
\end{equation}
and a fortiori \eqref{key-lem101} holds.
 
We then prove by induction on $p$ that, for every integer $p\geq 2$ and every
$\tau\in\mathbf{T}^b_p$,
\begin{equation}
\label{key-lem104}
100\,c_1^p\,\ve^{1+4p}(x\vee \ve)^{-7}\leq \int \Lambda^{(\tau)}(\dd \theta)\,J^\ve_x(\theta)\leq \ve^{1+4p}(x\vee \ve)^{-7}.
\end{equation}
So fix an integer $q\geq 2$, and 
 assume that \eqref{key-lem104} holds for $p=2,\ldots,q-1$ (when $q=2$ we make no assumption but we will
 rely on \eqref{key-lem103}). 
 Let $\tau\in\mathbf{T}^b_q$, and write $\tau_1$ and $\tau_2$ for the two binary trees obtained 
 from $\tau$ by removing the root. Then,
 \begin{align}
 \label{key-lem1}
 \int \Lambda^{(\tau)}(\dd \theta)\,J_x^\ve(\theta)
&=\int_0^\infty \dd t\,\E_x\Big[R_t^4\,\int \Lambda^{(\tau_1)}(\dd {\theta_1})\,J^\ve_{R_t}(\theta_1)\,\int \Lambda^{(\tau_2)}(\dd {\theta_2})\,J^\ve_{R_t}(\theta_2)\Big]\nonumber\\
&\leq \ve^{2+4q}\int_0^\infty \dd t\,\E_x\Big[ (R_t\vee \ve)^{-10}\Big],
 \end{align}
 where we applied either \eqref{key-lem103} or the induction hypothesis to both $\tau_1$ and $\tau_2$, using the trivial bound $(x\vee\ve)^{-7}\leq x^{-2}(x\vee\ve)^{-5}$ and noting that $\# L(\tau_1)+\# L(\tau_2)=\# L(\tau)= q$.

We then evaluate
 $$\int_0^\infty \dd t\,\E_x\Big[ (R_t\vee \ve)^{-10}\Big]
 = \int_{\R^9} \dd z\,(|z|\vee\ve)^{-10}\,G(z_x,z)
 =\frac{2}{7}\,\int_0^\infty\dd\rho\,\rho^{8}(\rho\vee \ve)^{-10}\,(x\vee\rho)^{-7},$$
 using the same argument as in the proof of \eqref{key-lem103}. Then, if $x\geq \ve$, we have
$$
\int_0^\infty\dd\rho\,\rho^{8}(\rho\vee \ve)^{-10}\,(x\vee\rho)^{-7}
 =\frac{1}{9}\ve^{-1}x^{-7} + (\ve^{-1}-x^{-1})x^{-7}+\frac{1}{8}x^{-8}
\leq \frac{10}{9} \ve^{-1}x^{-7}=\frac{10}{9} \ve^{-1}(x\vee \ve)^{-7},
$$
and, if $x<\ve$,
$$\int_0^\infty\dd\rho\,\rho^{8}(\rho\vee \ve)^{-10}\,(x\vee\rho)^{-7}
=\frac{1}{9}\ve^{-10}x^2+\frac{1}{2}(\ve^2-x^2)\ve^{-10}+\frac{1}{8}\ve^{-8}\leq \ve^{-8}= \ve^{-1}(x\vee \ve)^{-7}.$$
From the last three displays, we get that the right-hand side of \eqref{key-lem1} is bounded above by
$$\frac{2}{7}\times \ve^{2+4q}\times \frac{10}{9} \ve^{-1}(x\vee \ve)^{-7}\leq \ve^{1+4q}(x\vee \ve)^{-7},$$
and we have obtained the upper bound of \eqref{key-lem104} for $p=q$. 

To get the lower bound, we argue similarly.
The upper bound \eqref{key-lem1} is now replaced by
$$ \int \Lambda^{(\tau)}(\dd \theta)\,J_x^\ve(\theta) \geq 100^2\,c_1^q\,
 \ve^{2+4q}\int_0^\infty \dd t\,\E_x\Big[ R_t^4\,(R_t\vee \ve)^{-14}\Big],
$$
using \eqref{key-lem103} and the induction hypothesis. The same elementary calculations show that
$$\int_0^\infty \dd t\,\E_x\Big[ R_t^4\,(R_t\vee \ve)^{-14}\Big]
= \frac{2}{7}\,\int_0^\infty\dd\rho\,\rho^{12}(\rho\vee \ve)^{-14}\,(x\vee\rho)^{-7}
\geq  \frac{1}{100}\,\ve^{-1}(x\vee \ve)^{-7}.$$
By combining the last two displays, we get the lower bound of \eqref{key-lem104} for $p=q$. This completes
the proof of the lemma and of Proposition \ref{key-estim}.
\endproof

We now state a variant of the upper bounds in Proposition \ref{key-estim}, where
we 
condition on $W_*=0$, instead of considering the event where $W_*>0$. This will be useful in our proof of Proposition \ref{moment-normalized} as we will need to
condition on the value of $W_*$.

\begin{proposition}
\label{moment-condi0}
For every integer $p\geq 1$, for every $\ve,x>0$,
$$\N_x\Big(\Big(\int_0^\sigma \dd s\,\mathbf{1}_{\{\wh W_s\leq \ve\}}\Big)^p \,\Big|\, W_*=0\Big) \leq (16C_2)^p\,p!\,\ve^{4p},$$
where $C_2$ is the constant in Proposition \ref{key-estim}.
\end{proposition}

\proof From Proposition \ref{spine-dec}, it is immediate 
to verify that the mapping
$$x\mapsto \N_x\Big(\Big(\int_0^\sigma \dd s\,\mathbf{1}_{\{\wh W_s\leq \ve\}}\Big)^p \,\Big|\, W_*=0\Big)$$
is nondecreasing. Hence, without loss of generality, we may assume that $x>\ve$. By Proposition \ref{key-estim}, we have
$$\N_{x+\ve}\Bigg( \Big(\int_0^\sigma \dd s\,\mathbf{1}_{\{\wh W_s\leq 2\ve\}}\Big)^p\,\mathbf{1}_{\{0<W_*<\ve\}}\Bigg)
 \leq C_2^p\,p!\,(x+\ve)^{-3}(2\ve)^{1+4p}.$$
 On the other hand, the left-hand side of the last display is bounded below by
 \begin{align*}
 &\N_{x+\ve}\Big( \Big(\int_0^\sigma \dd s\,\mathbf{1}_{\{\wh W_s-W_*<\ve\}}\Big)^p\,\mathbf{1}_{\{0<W_*<\ve\}}\Big)\\
 &\quad = \int_{x}^{x+\ve} \dd u\,\frac{3}{u^3}
 \N_{x+\ve}\Big( \Big(\int_0^\sigma \dd s\,\mathbf{1}_{\{\wh W_s-W_*<\ve\}}\Big)^p\,\Big|\,W_*=x+\ve-u\Big)\\
 &\quad = \int_{x}^{x+\ve} \dd u\,\frac{3}{u^3}
 \N_{u}\Big( \Big(\int_0^\sigma \dd s\,\mathbf{1}_{\{\wh W_s<\ve\}}\Big)^p\,\Big|\,W_*=0\Big)\\
 &\quad \geq \frac{3}{2}\Big(x^{-2}-(x+\ve)^{-2}\Big)\,
  \N_{x}\Big( \Big(\int_0^\sigma \dd s\,\mathbf{1}_{\{\wh W_s<\ve\}}\Big)^p\,\Big|\,W_*=0\Big),
  \end{align*}
  where we used \eqref{min-snake} to condition with respect to $W_*$ and in the last line we also applied
 the previously mentioned monotonicity property. We have thus obtained that
  $$ \N_{x}\Big( \Big(\int_0^\sigma \dd s\,\mathbf{1}_{\{\wh W_s<\ve\}}\Big)^p\,\Big|\,W_*=0\Big)
  \leq \Big(\frac{3}{2}\Big(x^{-2}-(x+\ve)^{-2}\Big)\Big)^{-1}\times C_2^p\,p!\,(x+\ve)^{-3}(2\ve)^{1+4p}.$$
  To complete the proof, just notice that
  $$x^{-2}-(x+\ve)^{-2}= \frac{2\ve x+\ve^2}{x^2(x+\ve)^2}\geq \frac{2\ve}{(x+\ve)^3}.\eqno{\square}$$
  
  We can now prove Proposition \ref{moment-normalized} which was stated at the beginning of this section.

\proof[Proof of Proposition \ref{moment-normalized}] By formula \eqref{decomp-N},
$$\N_0\Big(\Big(\int_0^\sigma \dd s\,\mathbf{1}_{\{\wh W_s - W_* \leq 2\ve\}}\Big)^p\,\mathbf{1}_{\{1<\sigma<2\}}\Big)
= \int_1^2 \frac{\dd u}{2\sqrt{2\pi u^3}}\,\N^{(u)}_0\Big( \Big(\int_0^u \dd s\,\mathbf{1}_{\{\wh W_s - W_* \leq2\ve\}}\Big)^p\Big).$$
Furthermore, for every $u>0$,
\begin{align*}
\N^{(u)}_0\Big( \Big(\int_0^u \dd s\,\mathbf{1}_{\{\wh W_s - W_* \leq2\ve\}}\Big)^p\Big)
%&=u^p \,\N^{(u)}_0\Big( \Big(\int_0^1 \dd s\,\mathbf{1}_{\{\wh W_{us} - W_* \leq2\ve\}}\Big)^p\Big)\\
&=u^p \,\N^{(u)}_0\Big( \Big(\int_0^1 \dd s\,\mathbf{1}_{\{u^{-1/4}\wh W_{us} - u^{-1/4}W_* \leq2u^{-1/4}\ve \}}\Big)^p\Big)\\
&=u^p\,\N^{(1)}_0\Big(\Big(\int_0^1 \dd s\,\mathbf{1}_{\{\wh W_{s} - W_* \leq2u^{-1/4}\ve \}}\Big)^p\Big),
\end{align*}
since the pushforward of $\N^{(u)}_0$ under the scaling operator $\Theta_{u^{-1/2}}$ is $\N^{(1)}_0$ 
(see the end of Section \ref{sec:Excu-mea}). It follows that
\begin{align}
\label{norma-1}
\N_0\Big(\Big(\int_0^\sigma \dd s\,\mathbf{1}_{\{\wh W_s - W_* \leq 2\ve\}}\Big)^p\,\mathbf{1}_{\{1<\sigma<2\}}\Big)
%&= \int_1^2 \frac{\dd u}{2\sqrt{2\pi u^3}}\,\N^{(u)}_0\Big( \Big(\int_0^u \dd s\,\mathbf{1}_{\{\wh W_s - W_* \leq2\ve\}}\Big)^p\Big)\nonumber\\
&= \int_1^2 \frac{\dd u}{2\sqrt{2\pi u^3}}\,u^p\,\N^{(1)}_0\Big( \Big(\int_0^1 \dd s\,\mathbf{1}_{\{\wh W_s - W_* \leq 2u^{-1/4}\ve \}}\Big)^p\Big)\nonumber\\
&\geq \frac{1}{4\sqrt{2\pi}}\,\N^{(1)}_0\Big( \Big(\int_0^1 \dd s\,\mathbf{1}_{\{\wh W_s - W_* \leq \ve\}}\Big)^p\Big).
\end{align}

In order to prove Proposition \ref{moment-normalized}, it is thus enough to bound the left-hand side of \eqref{norma-1}. 
By conditioning with respect to $W_*$, using \eqref{min-snake}, we get
 $$\N_0\Big(\Big(\int_0^\sigma \dd s\,\mathbf{1}_{\{\wh W_s - W_* \leq \ve\}}\Big)^p\,\mathbf{1}_{\{1<\sigma<2\}}\Big)
 =3\int_0^\infty \frac{\dd x}{x^3}\,\N_0\Big(\Big(\int_0^\sigma \dd s\,\mathbf{1}_{\{\wh W_s - W_* \leq \ve\}}\Big)^p\,\mathbf{1}_{\{1<\sigma<2\}}
 \,\Big|\,W_*=-x\Big).$$
 When $x\geq 1$, we use Proposition \ref{moment-condi0} to bound
 \begin{align}
 \label{norma-2}
 \N_0\Big(\Big(\int_0^\sigma \dd s\,\mathbf{1}_{\{\wh W_s - W_* \leq \ve\}}\Big)^p\,\mathbf{1}_{\{1<\sigma<2\}}
 \,\Big|\,W_*=-x\Big)
 &\leq \N_0\Big(\Big(\int_0^\sigma \dd s\,\mathbf{1}_{\{\wh W_s - W_* \leq \ve\}}\Big)^p
 \,\Big|\,W_*=-x\Big)\nonumber\\
 &=\N_x\Big(\Big(\int_0^\sigma \dd s\,\mathbf{1}_{\{\wh W_s\leq \ve\}}\Big)^p
 \,\Big|\,W_*=0\Big)\nonumber\\
 &\leq (16C_2)^p\,p!\,\ve^{4p}.
 \end{align}
 Suppose then that $0<x<1$. By the H\"older inequality,
 \begin{align}
 \label{norma-3}
 &\N_0\Big(\Big(\int_0^\sigma \dd s\,\mathbf{1}_{\{\wh W_s - W_* \leq \ve\}}\Big)^p\,\mathbf{1}_{\{1<\sigma<2\}}
 \,\Big|\,W_*=-x\Big)\nonumber\\
&\qquad \leq \Bigg(\N_0\Big(\Big(\int_0^\sigma \dd s\,\mathbf{1}_{\{\wh W_s - W_* \leq \ve\}}\Big)^{4p}
 \,\Big|\,W_*=-x\Big)\Bigg)^{1/4} \times \Big(\N_0(1<\sigma<2\mid W_*=-x)\Big)^{3/4}.
 \end{align}
Using again Proposition \ref{moment-condi0}, we get the existence of a 
constant $C_3$ such that
\begin{equation}
\label{norma-4}
\Bigg(\N_0\Big(\Big(\int_0^\sigma \dd s\,\mathbf{1}_{\{\wh W_s - W_* \leq \ve\}}\Big)^{4p}
 \,\Big|\,W_*=-x\Big)\Bigg)^{1/4}
 \leq \Big( (16C_2)^{4p} \,(4p)!\,\ve^{16p}\Big)^{1/4}
 \leq C_3^p\,p!\,\ve^{4p}.
\end{equation}
On the other hand, we have
$$\N_0(1<\sigma<2\mid W_*=-x)\leq \frac{1}{1-e^{-1/2}}\,\N_0(1-e^{-\sigma/2}\mid W_*=-x),$$
and by \cite[Proposition 4.6]{CLG},
$$\N_0(1-e^{-\sigma/2}\mid W_*=-x)= 1- \frac{x^3\cosh(x)}{(\sinh(x))^3}= O(x^4)$$
as $x\to 0$. Hence there exists a constant $c'$ such that, for every $x\in(0,1)$,
\begin{equation}
\label{norma-5}
\N_0(1<\sigma<2\mid W_*=-x)\leq c'\,x^4.
\end{equation}
Finally, using \eqref{norma-2},\eqref{norma-3},\eqref{norma-4},\eqref{norma-5}, we obtain that
 $$\N_0\Big(\Big(\int_0^\sigma \dd s\,\mathbf{1}_{\{\wh W_s - W_* \leq \ve\}}\Big)^p\,\mathbf{1}_{\{1<\sigma<2\}}\Big)
 \leq 3\int_1^\infty \frac{\dd x}{x^3}\,(16C_2)^p\,p!\,\ve^{4p} + 3\int_0^1 \frac{\dd x}{x^3}\,C_3^p\,p!\,\ve^{4p}\, (c'\,x^4)^{3/4},$$
 and the bound of Proposition \ref{moment-normalized} follows from \eqref{norma-1}. \endproof
 
 Recall that $\bm_\infty$ stands for the (standard) Brownian sphere defined
 under the probability measure $\N^{(1)}_0$, and that $\mathrm{Vol}(\cdot)$ denotes the volume 
 measure on $\bm_\infty$. Also recall that, for $a\in\bm_\infty$ and $r>0$, $B(a,r)$ denotes the closed ball of radius $r$ centered 
 at $a$ in $\bm_\infty$.
 
 \begin{proposition}
 \label{lower-Haus}
 For every $r\in(0,1/4]$, set
 $h(r)=r^4 \log\log(1/r)$. There exists a constant $K_1$ such that we have $\N^{(1)}_0$ a.s.
 $$\mathrm{Vol}(\dd a)\hbox{ a.e.}, \quad\limsup_{r\downarrow  0} \frac{\mathrm{Vol}(B(a,r))}{h(r)} \leq K_1.$$
 Consequently, there exists a constant $\kappa_1>0$ such that we have $\N^{(1)}_0$ a.s. for every Borel subset $A$ of $\bm_\infty$, 
 $$m^h(A)\geq \kappa_1\,\mathrm{Vol}(A).$$
 \end{proposition}
 
 \proof The first assertion will follow from
Proposition \ref{re-rooting} if we can verify that, for some constant $K_1$,
 \begin{equation}
 \label{lower-Haus1}
 \limsup_{r\downarrow  0} \frac{\mathrm{Vol}(B(\mathbf{x}_*,r))}{h(r)} \leq K_1,\quad \hbox{a.s.}
 \end{equation}
 Let $\lambda\in(0,1/C_0)$. By \eqref{vol-ball-x*} and Proposition \ref{moment-normalized}, we have
 for every $r>0$,
 $$\N^{(1)}_0\Big(\exp\Big(\lambda \frac{\mathrm{Vol}(B(\mathbf{x}_*,r)}{r^4}\Big)\Big) \leq K_{(\lambda)},$$
 where $K_{(\lambda)}$ is a finite constant depending on $\lambda$. Hence, for every $r>0$ and $u\geq 0$,
 we have by the Markov inequality,
 $$\N^{(1)}_0(\mathrm{Vol}(B(\mathbf{x}_*,r))>u\,r^4)\leq K_{(\lambda)}\,\exp(-\lambda u).$$
Apply this bound with $\lambda= 1/(2C_0)$, $r=2^{-k}$ (for $k\in \N$) and 
$u=4C_0\log k$. The Borel-Cantelli lemma then gives that a.s. for $k$ large enough we have
$$\mathrm{Vol}(B(\mathbf{x}_*,2^{-k})) \leq 4C_0\,2^{-4k}\,\log k,$$
and \eqref{lower-Haus1} follows. 

To get the second assertion, if $A$ is a Borel subset of $\bm_\infty$, we set
$$A':=\Bigg\{a\in A: \limsup_{r\downarrow  0} \frac{\mathrm{Vol}(B(a,r))}{h(r)} \leq K_1\Bigg\}$$
so that $\mathrm{Vol}(A')=\mathrm{Vol}(A)$ (by the first assertion) and $m^h(A')\leq m^h(A)$.
Applying Lemma \ref{compa}(i) to $A'$ yields the desired result.
\endproof

\section{The upper bound for the Hausdorff measure}
\label{sec:upper}

Proposition \ref{lower-Haus} provides a lower bound for the Hausdorff measure of a 
Borel subset $A$ of $\bm_\infty$. Our goal in this section is to derive the
(more delicate) upper bound. More precisely, we will prove the following 
proposition.

\begin{proposition}
\label{upper-Haus}
Let $h$ be as in Proposition \ref{lower-Haus}.
There exists a constant $K_2>0$ such that we have $\N^{(1)}_0$ a.s.
$$m^h\Bigg(\Bigg\{a\in\bm_\infty: \limsup_{r\downarrow 0}\frac{\mathrm{Vol}(B(a,r))}{h(r)} < K_2\Bigg\}\Bigg) =0.$$
Consequently, there exists a constant $\kappa_2>0$ such that, $\N^{(1)}_0$ a.s. for every Borel subset $A$ of $\bm_\infty$, 
 $$m^h(A)\leq \kappa_2\,\mathrm{Vol}(A).$$
\end{proposition}

As in the case of Proposition \ref{lower-Haus}, the second part of the proposition almost immediately follows 
from the first assertion. So the difficult part is to prove the first assertion, and, to this end, we will rely on
several lemmas. We need to control the Hausdorff measure of the set of all points $a$ of $\bm_\infty$
such that the volumes of balls centered at $a$ are unusually small. In a way similar to the previous section, 
we start by considering balls centered at the distinguished point $\bx_*$, for which we can use formula 
\eqref{vol-ball-x*}. The key estimate needed to handle these balls is Lemma \ref{upper-key2}, which bounds 
the probability that the volumes of the balls $B(\bx_*,2^{-k})$ are small for all $\lfloor m/2\rfloor\leq k\leq m$. As in the previous section,
although Lemma \ref{upper-key2} is stated under the probability measure $\N^{(1)}_0$, it is more convenient to
start by proving similar estimates under the measure $\N_x(\cdot\mid W_*=0)$, and this is the motivation 
for Lemmas \ref{upper-lem1} and \ref{upper-key}. Finally, in order to derive Proposition \ref{upper-Haus}
from Lemma \ref{upper-key2}, we rely both on the re-rooting invariance property (Proposition \ref{re-rooting})
and on a uniform modulus of continuity for the pseudo-distance $D(s,t)$. We start by stating and proving the 
latter estimate. 

\begin{lemma}
\label{Holder}
$\N^{(1)}_0(\dd\omega)$ a.s., there exists a finite constant $C(\omega)$ such that, for every distinct $s,t\in[0,1]$,
$$|\wh W_s-\wh W_t|\leq C(\omega)\,\Big(1+ \log\frac{1}{|t-s|}\Big)\,|t-s|^{1/4},$$
and
$$D(s,t)\leq C(\omega)\,\Big(1+ \log\frac{1}{|t-s|}\Big)\,|t-s|^{1/4}.$$
\end{lemma}

\proof It is enough to prove the first bound of the lemma. Indeed,
$$D(s,t)\leq D^\circ(s,t)\leq \wh W_s + \wh W_t -2\,\min_{s\wedge t\leq u\leq s\vee t} \wh W_u,$$
and the right-hand side can be written as
$|\wh W_s-\wh W_u|+|\wh W_t-\wh W_u|$ for some $u\in[s\wedge t,s\vee t]$. 

In view of deriving the first bound, we start by considering the increments
$|\wh W'_s-\wh W'_t|$, when $(W'_u)_{u\geq 0}$ is a (one-dimensional) Brownian snake starting from $0$ under the 
probability measure $\P$ (see e.g.~\cite{Zurich}). Writing $\beta_u=\zeta_{(W'_u)}$ to simplify notation, this means in particular that
the process $(\beta_u)_{u\geq 0}$
is a reflecting Brownian motion on the half-line, and that, conditionally on $(\beta_u)_{u\geq 0}$,
$(\wh W'_u)_{u\geq 0}$ is a centered Gaussian process with continuous sample paths and covariance
$$\E[\wh W'_s\wh W'_t\mid (\beta_u)_{u\geq 0}]= \min_{s\wedge t\leq u\leq s\vee t} \beta_u.$$
Fix $s,t\geq 0$ with $s\not =t$ and set
$$d_\beta(s,t):=\beta_s+\beta_t- 2\min_{s\wedge t\leq u\leq s\vee t} \beta_u.$$
Then, conditionally on $(\beta_u)_{u\geq 0}$, $\wh W'_s-\wh W'_t$
is a centered Gaussian variable with variance $d_\beta(s,t)$, and thus, for every $\lambda>0$,
$$\E\Big[ \exp(\lambda(\wh W'_s-\wh W'_t))\,\Big|\, (\beta_u)_{u\geq 0}\Big]= \exp(\frac{\lambda^2}{2}\,d_\beta(s,t)).$$
Furthermore, we leave it as an exercise for the reader to check that $d_\beta(s,t)$ is stochastically dominated
by $R^{(3)}_{|t-s|}$, where $(R^{(3)}_u)_{u\geq 0}$ stands for a three-dimensional Bessel
process started from $0$ (if $\beta$ were a linear Brownian motion, the celebrated Pitman theorem would 
show that the law of $d_\beta(s,t)$ is exactly the law of $R^{(3)}_{|t-s|}$). We thus get, for every
$\lambda>0$,
$$\E\Big[ \exp(\lambda(\wh W'_s-\wh W'_t))\Big]=\E[\exp(\frac{\lambda^2}{2}\,d_\beta(s,t))]\leq \E\Big[\exp(\frac{\lambda^2}{2} R^{(3)}_{|t-s|})\Big]
= \E\Big[\exp(\frac{\lambda^2}{2}\,|t-s|^{1/2}\, R^{(3)}_{1})\Big],
$$
by scaling. We take $\lambda=|t-s|^{-1/4}$ and use the Markov inequality to get
for every $\alpha>0$,
$$\P(\wh W'_s-\wh W'_t > \alpha |t-s|^{1/4}) \leq e^{-\alpha} \,\E\Big[\exp\Big(|t-s|^{-1/4}(\wh W'_s-\wh W'_t)\Big)\Big]
\leq e^{-\alpha} \E[\exp(R^{(3)}_1/2)]= \ov{c}\ e^{-\alpha},$$
where $\ov{c}= \E[\exp(R^{(3)}_1/2)]$ is a constant. 

It follows from the last bound that, for every integer $p\geq 1$ and $j\in\{1,2,\ldots,2^p\}$, we have
$$\P\Big(|\wh W'_{j2^{-p}}-\wh W'_{(j-1)2^{-p}}|>p2^{-p/4}\Big)\leq 2\ov{c}\,e^{-p}.$$
Summing over $j=1,\ldots,2^p$ and using the Borel-Cantelli lemma, we get that a.s. there exists
an integer $p_0$ such that, for every $p\geq p_0$ and $j\in\{1,2,\ldots,2^p\}$, we have
$$|\wh W'_{j2^{-p}}-\wh W'_{(j-1)2^{-p}}|\leq p2^{-p/4}.$$
By  standard chaining arguments similar to the proof
of the classical Kolmogorov lemma, we derive from the latter bound that 
$$\sup_{s,t\in[0,1], s\not = t} \frac{|\wh W'_s-\wh W'_t|}{(1+\log(1/|t-s|))|t-s|^{1/4}} <\infty, \quad \P\hbox{ a.s.}$$
Since $\N_0$ can be interpreted as the excursion measure away from $0$ for the Brownian snake, the
result of the lemma follows from the last display.
\endproof

\begin{lemma}
\label{upper-lem1}
There exist two positive constants $\wt c$ and $\alpha_0$ such that, for every $\ve>0$ and $x\in[2\ve,3\ve]$,
we have for every $u\geq 0$,
$$\N_x\Big(\int_0^\sigma \dd s\,\mathbf{1}_{\{\wh W_s\leq \ve\}}>u\ve^4\,\Big|\, 0<W_*\leq \ve\Big)\geq \wt c\, e^{-\alpha_0u}.$$
\end{lemma}

\proof This is basically a consequence of Proposition \ref{key-estim}. From this proposition 
we get, for every $x\in[2\ve,3\ve]$, for every integer $p\geq 1$,
$$\frac{1}{27}C_1^pp!\,\ve^{4p-2}\leq\N_x\Bigg( \Big(\int_0^\sigma \dd s\,\mathbf{1}_{\{\wh W_s\leq \ve\}}\Big)^p\,\mathbf{1}_{\{0<W_*\leq\ve\}}\Bigg)\leq 
\frac{1}{8}C_2^p p!\,\ve^{4p-2}.$$
Then note that, for $x\in[2\ve,3\ve]$, the quantity
$$\N_x(0<W_*\leq\ve)=\frac{3}{2}((x-\ve)^{-2}-x^{-2})$$
is bounded above by $9\ve^{-2}/8$ and bounded below by $5\ve^{-2}/24$. Hence, we get the existence of (explicit) constants
$\wt c_1$ and $\wt c_2$ such that, for every $\ve>0$ and $x\in[2\ve,3\ve]$, we have for every integer $p\geq1$,
$$\wt c_1\,C_1^p\,p!\leq \N_x\Bigg( \Big(\ve^{-4}\int_0^\sigma \dd s\,\mathbf{1}_{\{\wh W_s\leq \ve\}}\Big)^p\,\Bigg|\, 0<W_*\leq\ve\Bigg)\leq \wt c_2\,C_2^p\,p!.$$
This estimate on moments implies the tail estimate of the lemma, with constants $\wt c$ and $\alpha_0$
that only depend on $\wt c_1,\wt c_2, C_1, C_2$. To see this, set $X_\ve=\ve^{-4}\int_0^\sigma \dd s\,\mathbf{1}_{\{\wh W_s\leq \ve\}}$
and also write $\P_{(\ve)}$, resp.~$\E_{(\ve)}$, for the probability measure $\N_x(\cdot\mid 0<W_*\leq \ve)$, resp.~for the expectation under
this probability measure. Then an application of the Cauchy-Schwarz inequality gives for every integer $p\geq 1$ and every $u>0$,
$$\E_{(\ve)}[X_\ve^p] \leq u^p+ \E_{(\ve)}[X_\ve^{2p}]^{1/2}\, \P_{(\ve)}(X_\ve >u)^{1/2}$$
and therefore
$$\P_{(\ve)}(X_\ve>u)^{1/2}\geq \frac{\big(\E_{(\ve)}[X_\ve^p]-u^p\big)^+}{\E_\ve[X_\ve^{2p}]^{1/2}}\geq \frac{\big(\wt c_1C_1^p p! - u^p\big)^+}
{\wt c_2C_2^{p}\,((2p)!)^{1/2}}.$$
Applying this bound with $p=\lceil 4u/C_1\rceil$ leads to the desired estimate. 
\endproof

In the remaining part of this section, we use the notation 
$$\v_\ve:=\int_0^\sigma \dd s\,\mathbf{1}_{\{\wh W_s-W_*\leq\ve\}}$$ 
for $\ve>0$. Under the probability measure $\N^{(1)}_0$, 
$\v_\ve$ is the volume of the ball of radius $\ve$ centered at the
distinguished point $\bx_*$ of the Brownian sphere (cf.~\eqref{vol-ball-x*}).

Recall the 
definition of the function $h$ in Proposition \ref{lower-Haus}.

\begin{lemma}
\label{upper-key}
There exist constants $\gamma>0$, $K>0$ and $\delta\in(0,1)$ such that, for every $x>0$
and every integer $n\geq 4$ such that $x>2^{-\lfloor n/2\rfloor}$, one has
$$\N_x\Bigg(\bigcap_{k=\lfloor n/2\rfloor}^n \{\v_{2^{-k}}\leq \gamma\,h(2^{-k})\}\,\Bigg|\, W_*=0\Bigg) \leq \exp(-K\,n^\delta).$$
\end{lemma}

\proof Let $x>0$. We rely on the spine decomposition under $\N_x(\cdot\midd W_*=0)$, which is given in Proposition \ref{spine-dec}.
To simplify notation, and only in this proof, we write $\P$ for the probability measure $\N_x(\cdot\midd W_*=0)$.
Let $U=(U_t)_{0\leq t\leq L_x}$ and 
$\mathcal{N}(\dd t\,\dd\omega')=\sum_{i\in I} \delta_{(t_i,\omega_i)}(\dd t\,\dd\omega')$
be defined under $\P$ as in Proposition \ref{spine-dec}. Then, for 
every $\ve>0$, we have $\P$ a.s.,
$$\v_\ve=\sum_{i\in I} \int_0^{\sigma(\omega_i)}\dd s\,\mathbf{1}_{\{\wh W_s(\omega_i)\leq \ve\}}.$$

Let $k_0$ be the smallest integer such that $k_0\geq 2$ and $2^{-k_0}<x$. Set
$L_{2^{-k}}:=\sup\{t\in[0,L_x]:U_t=2^{-k}\}$ for every $k\geq k_0$, and $L^{(k)}:=L_{2^{-k}}-L_{2^{-k-1}}$. Then the processes $U^{(k)}$ defined for $k\geq k_0$ by
$$U^{(k)}_t:=U_{(L_{2^{-k-1}}+t)\wedge L_{2^{-k}}},\quad t\geq 0,$$
are independent, and furthermore the distribution of $(2^kU^{(k)}_{2^{-2k}t})_{0\leq t\leq 2^{2k}L^{(k)}}$ does not
depend on $k$. The preceding facts are consequences of well-known properties 
of Bessel processes. It follows that the point measures
$$\mathcal{N}_k(\dd t\,\dd\omega'):=\sum_{i\in I,\,L_{2^{-k-1}}\leq t_i<L_{2^{-k}}} \delta_{(t_i-L_{2^{-k-1}},\omega_i)}(\dd t\, \dd\omega')$$
are independent when $k$ varies in $\{k_0,k_0+1,\ldots\}$. 

 For every $k\geq k_0$, let $A_k$ denote the event where the point measure $\mathcal{N}_k$ has an atom $(\ov t,\ov\omega)$ such 
 that $\ov\omega(0)\in [2\times 2^{-k},3\times 2^{-k}]$
 and $W_*(\ov\omega)<2^{-k}$. On the event $A_k$, write $(t_{(k)},\omega_{(k)})$ for this atom
 (if there are several possible choices, choose the one with the smallest value 
 of $\ov t$). By the last property stated in Proposition \ref{spine-dec}, we know that, conditionally on $U$, $\mathcal{N}_k(\dd t\,\dd\omega')$ is Poisson on
 $[0,L^{(k)}] \times \S$ 
 with intensity 
 $$4\,\dd t\,\N_{U^{(k)}_t}
(\dd \omega' \cap \{W_*(\omega')>0\}).$$
Hence, the probability of the complement of $A_k$ is 
\begin{align*}
\P\big((A_k)^c\big)&=\E\Bigg[ \exp\Bigg(-4\int_{0}^{L^{(k)}} \dd t\,\mathbf{1}_{\{U^{(k)}_t\in[2\times 2^{-k},3\times 2^{-k}]\}}
 \times \N_{U^{(k)}_t}\Big(0<W_*<2^{-k}\Big)\Bigg)\Bigg]\\
&= \E\Bigg[ \exp\Bigg(-4\int_0^{L^{(k)}} \dd t\,\mathbf{1}_{\{U^{(k)}_t\in[2\times 2^{-k},3\times 2^{-k}]\}}
 \times \frac{3}{2}\Big((U^{(k)}_t-2^{-k})^{-2}-(U^{(k)}_t)^{-2}\Big)\Bigg)\Bigg]\\
 &= \E\Bigg[ \exp\Bigg(-4\,2^{2k}\int_0^{L^{(k)}} \dd t\,\mathbf{1}_{\{2^{k}U^{(k)}_t\in[2,3]\}}
 \times \frac{3}{2}\Big((2^kU^{(k)}_t-1)^{-2}-(2^kU^{(k)}_t)^{-2}\Big)\Bigg)\Bigg]
 \end{align*}
 and the change of variable $t=2^{-2k}u$ (together with the fact that the distribution of $(2^kU^{(k)}_{2^{-2k}t})_{0\leq t\leq 2^{2k}L^{(k)}}$ does not
depend on $k$) shows that the right-hand side of the last display does not depend on $k$. We set
 $\eta:=\P(A_k)$, which is a positive constant independent of $k$. We also observe that,
under $\P(\cdot\,|\,A_k)$ and conditionally on $\omega_{(k)}(0)=y\in [2\times 2^{-k},3\times 2^{-k}]$, $\omega_{(k)}$
 is distributed according to $\N_y(\cdot\mid 0<W_*\leq 2^{-k})$, as a consequence of properties of
 Poisson measures. 
 
 Let $n\geq 4$ be an integer such that $x>2^{-\lfloor n/2\rfloor}$, so that $k_0\leq \lfloor n/2\rfloor$. Since we have trivially, on the event $A_k$,
$$\mathcal{V}_{2^{-k}}\geq \int_0^{\sigma(\omega_{(k)})}\dd s
 \,\mathbf{1}_{\{\wh W_s(\omega_{(k)})\leq 2^{-k}\} },$$
we get for every $\gamma>0$,
 \begin{align*}
 \P\Bigg(\bigcap_{k=\lfloor n/2\rfloor}^n \{\v_{2^{-k}}\leq \gamma\,h(2^{-k})\}\Bigg)
 &\leq  \P\Bigg(\bigcap_{k=\lfloor n/2\rfloor}^n\Big(A_k\cap \Big\{\int_0^{\sigma(\omega_{(k)})}\dd s
 \,\mathbf{1}_{\{\wh W_s(\omega_{(k)})\leq 2^{-k}\} }>\gamma\,h(2^{-k})\Big\}\Big)^c\Bigg)\\
 &= \prod_{k=\lfloor n/2\rfloor}^n \Bigg(1- \P\Big(A_k\cap \Big\{\int_0^{\sigma(\omega_{(k)})}\dd s
 \,\mathbf{1}_{\{\wh W_s(\omega_{(k)})\leq 2^{-k}\}} >\gamma\,h(2^{-k})\Big\}\Big)\Bigg),
 \end{align*}
thanks to the independence of the point measures $\mathcal{N}_k$. Then, using the conditional
distribution 
of $\omega_{(k)}$, we have for every $k\geq k_0$,
\begin{align*}
&\P\Big(A_k\cap \Big\{\int_0^{\sigma(\omega_{(k)})}\dd s
 \,\mathbf{1}_{\{\wh W_s(\omega_{(k)})\leq 2^{-k}\}} >\gamma\,h(2^{-k})\Big\}\Big)\\
&\quad \geq \eta \;\inf_{y\in[2.2^{-k},3.2^{-k}]}
 \N_y\Big( \int_0^{\sigma}\dd s\,\mathbf{1}_{\{\wh W_s\leq 2^{-k}\}}>\gamma\,h(2^{-k})\,\Big|\,0<W_*\leq 2^{-k}\Big)
 \\
 &\quad\geq \eta\,\times \wt c\,\exp(-\alpha_0\,\gamma\,\log\log 2^k)
 \end{align*}
by Lemma \ref{upper-lem1}. Now fix $\gamma$ such that $\alpha_0\,\gamma<1$. The right-hand side of the 
last display is bounded below by $\tilde c'\,k^{-\alpha_0\gamma}$ for some constant $\tilde c'>0$. It follows that
$$
 \P\Bigg(\bigcap_{k=\lfloor n/2\rfloor}^n \{\v_{2^{-k}}\leq \gamma\,h(2^{-k})\}\Bigg)
 \leq \prod_{k=\lfloor n/2\rfloor}^n \Big(1-\tilde c'\,k^{-\alpha_0\gamma}\Big)
 \leq \exp\Big(-K\,n^{1-\alpha_0\gamma}\Big)
$$
 with some positive constant $K$. This completes the proof of the lemma. \endproof
 
 We can now use Lemma \ref{upper-key} to derive the estimate that we need under $\N^{(1)}_0$.
 
 \begin{lemma}
 \label{upper-key2}
 Let $\gamma$, $K$ and $\delta$ be as in Lemma \ref{upper-key}. There exists a constant $K'$
 such that, for every integer $n\geq 4$,
 $$\N^{(1)}_0\Bigg(\bigcap_{k=\lfloor n/2\rfloor}^n \{\v_{2^{-k}}\leq \frac{\gamma}{2}\,h(2^{-k})\}\Bigg) \leq K'\exp(-\frac{K}{4}\,n^\delta).$$
\end{lemma}

\proof To simplify notation, we write $E_{n,\gamma}$ for the event
$$E_{n,\gamma}:=\bigcap_{k=\lfloor n/2\rfloor}^n \{\v_{2^{-k}}\leq \gamma\,h(2^{-k})\}.$$
We have
\begin{equation}
\label{upper-tech10}
\N_0(E_{n,\gamma}\cap \{1<\sigma<2\}) = \int_1^2 \frac{\dd r}{2\sqrt{2\pi r^3}}\,\N_0^{(r)}(E_{n,\gamma}).
\end{equation}
We then use a scaling argument. Recall the definition of the scaling operators $\Theta_\lambda$
at the end of Section \ref{sec:Excu-mea}, and the fact that $\Theta_{\lambda}(\omega)$ is distributed
according to $\N^{(\lambda^2)}_0$ if $\omega$ is distributed according to $\N^{(1)}_0$.
For $r>0$, write $\omega^{(r)}=\Theta_{\sqrt{r}}(\omega)$ to simplify notation.
Then, for every $\ve>0$,
$$\v_\ve(\omega^{(r)})=\int_0^r \dd s\,\mathbf{1}_{\{\wh\omega^{(r)}_s-\omega^{(r)}_*\leq \ve\}}
= r \int_0^1 \dd s'\,\mathbf{1}_{\{\wh\omega_{s'}-\omega_*\leq \ve\,r^{-1/4}\}}=r\v_{\ve r^{-1/4}}(\omega).$$
It follows that, for $1<r<2$,
\begin{align*}
\N_0^{(r)}(E_{n,\gamma})=\N_0^{(r)}\Bigg(\bigcap_{k=\lfloor n/2\rfloor}^n \{\v_{2^{-k}}\leq \gamma\,h(2^{-k})\}\Bigg)
&=\N_0^{(1)}\Bigg(\bigcap_{k=\lfloor n/2\rfloor}^n \{\v_{r^{-1/4}2^{-k}}\leq \frac{\gamma}{r}\,h(2^{-k})\}\Bigg)\\
&\geq \N_0^{(1)}\Bigg(\bigcap_{k=\lfloor n/2\rfloor}^n \{\v_{2^{-k}}\leq \frac{\gamma}{2}\,h(2^{-k})\}\Bigg)\\
&= \N_0^{(1)}(E_{n,\gamma/2}),
\end{align*}
using the simple fact $\v_{r^{-1/4}2^{-k}}\leq \v_{2^{-k}}$ if $r>1$.
By substituting the latter bound in \eqref{upper-tech10}, we get
\begin{equation}
\label{upper-tech11}
\N_0^{(1)}(E_{n,\gamma/2})\leq \Big(\int_1^2 \frac{\dd r}{2\sqrt{2\pi r^3}}\Big)^{-1}\,\N_0(E_{n,\gamma}\cap \{1<\sigma<2\}).
\end{equation}
It remains to bound $\N_0(E_{n,\gamma}\cap \{1<\sigma<2\})$. Conditioning with respect to $W_*$, we have
\begin{equation}
\label{upper-tech12}
\N_0(E_{n,\gamma}\cap \{1<\sigma<2\})= 3\int_0^\infty \frac{\dd x}{x^3}\,\N_0(E_{n,\gamma}\cap \{1<\sigma<2\}\mid W_*=-x).
\end{equation}
For $x>2^{-\lfloor m/2\rfloor}$, Lemma \ref{upper-key} gives
$$\N_0(E_{n,\gamma}\mid W_*=-x)\leq \exp(-K\,n^\delta),$$
and on the other hand, we have also by \eqref{norma-5}, for $0<x<1$,
$$\N_0(1<\sigma<2\mid W_*=-x)\leq c'\,x^4.$$
From these estimates and \eqref{upper-tech12}, we get 
\begin{align*}
\frac{1}{3}\N_0(E_{n,\gamma}\cap \{1<\sigma<2\})&\leq c'\int_0^{2^{-\lfloor n\rfloor/2}} x\,\dd x
+ \int_{2^{-\lfloor n\rfloor/2}}^1 (\exp(-Kn^\delta)\wedge c'x^4)\,\frac{\dd x}{x^3} + \exp(-Kn^\delta)\int_1^\infty \frac{\dd x}{x^3}\\
&\leq c'\,2^{-n} + (c')^{3/4}\,\exp(-\frac{K}{4}n^\delta) + \frac{1}{2}\exp(-Kn^\delta),
\end{align*}
using the trivial bound $\exp(-Kn^\delta)\wedge c'x^4\leq (c')^{3/4}\,x^3\exp(-\frac{K}{4}n^\delta)$. The bound of the
lemma now follows from \eqref{upper-tech11} and the last display. \endproof

\noindent{\it Proof of Proposition \ref{upper-Haus}}. We prove that the first assertion of Proposition \ref{upper-Haus}
holds with $K_2=\gamma/32$, where $\gamma$ is as in the preceding two lemmas. To this end, it is 
enough to verify that, for every integer $n_0\geq 2$, we have $m^h(\b_{n_0})=0$ a.s., where
$$\b_{n_0}:=\Big\{ a\in\bm_\infty: \mathrm{Vol}(B(a,2^{-k}))\leq \frac{\gamma}{32}\,h(2^{-k}),\,\forall k\geq n_0\Big\}.$$
Let us fix the positive integer $n_0\geq 2$, and, for every $n>n_0$, set
$$\wt\b_{n_0,n}:=\Big\{ a\in\bm_\infty: \mathrm{Vol}(B(a,2^{-k}))\leq \frac{\gamma}{2}\,h(2^{-k}),\,\forall k\in\{n_0,n_0+1,\ldots,n\}\Big\}.$$
Recall the notation $\bp$ for the canonical projection
from $[0,1]$ onto $\bm_\infty$. By Proposition \ref{re-rooting}, for every integer $p\geq 1$ and every $i\in\{1,\ldots,2^p\}$,
$$\N^{(1)}_0\Big(\bp(i2^{-p})\in \wt\b_{n_0,n}\Big)=\N^{(1)}_0\Big(\mathbf{x}_*\in \wt\b_{n_0,n}\Big),$$
Suppose that $n\geq 2n_0$. Since $\mathrm{Vol}(B(\mathbf{x}_*,\ve))=\v_\ve$ for every $\ve >0$,
Lemma \ref{upper-key2} gives
$$\N^{(1)}_0\Big(\mathbf{x}_*\in \wt\b_{n_0,n}\Big)\leq K'\exp(-\frac{K}{4}\,n^\delta).$$
Set $\Gamma_{n_0,n}^p:=\{i\in\{1,2,\ldots,2^{p}\}: \bp(i2^{-p})\in \wt\b_{n_0,n}\}$. We get
$$\N^{(1)}_0(\#\Gamma_{n_0,n}^p) \leq K'2^p\,\exp(-\frac{K}{4}\,n^\delta).$$
Set $K''=K/8$. We apply the preceding bound with $n=\lfloor p/8\rfloor$
(when $p$ is large so that $p/8>2n_0$, which we assume
from now on). It immediately follows that
\begin{equation}
\label{upper-Haus2}
2^{-p}\,\exp\Big(K''(\frac{p}{8})^\delta\Big)\,\#\Gamma_{n_0,\lfloor p/8\rfloor}^p \build{\la}_{p\to\infty}^{} 0, \quad \N^{(1)}_0\hbox{ a.s.}
\end{equation}

The point is now to observe that, if $t\in[(i-1)2^{-p},i2^{-p}]$, for some $i\in\{1,\ldots,2^p\}$, we have
$$B(\bp(t),2^{-k})\supset B(\bp(i2^{-p}), 2^{-k-1})$$
as soon as $k$ is such that $D(i2^{-p},t)<2^{-k-1}$, which holds by Lemma \ref{Holder}
if 
$$C(\omega)\,\Big(1+\log 2^p\Big)\,2^{-p/4}<2^{-k-1},$$
and in particular if $k\leq \lfloor\frac{p}{8}\rfloor$, provided $p$ is large enough
(depending on $\omega$). From the preceding discussion we get that, for $p$
large enough, if $a\in \bm_\infty$ is of the form $a=\bp(t)$
with $t\in[(i-1)2^{-p},i2^{-p}]$, the condition $a\in \b_{n_0}$ implies that
for $n_0\leq k\leq \lfloor\frac{p}{8}\rfloor$,
$$\mathrm{Vol}\Big(B(\bp(i2^{-p}), 2^{-k-1})\Big)\leq \mathrm{Vol}\Big(B(a,2^{-k})\Big)
\leq \frac{\gamma}{32}\,h(2^{-k})\leq \frac{\gamma}{2}h(2^{-k-1}),$$
and therefore $\bp(i2^{-p})\in \wt \b_{n_0+1,\lfloor p/8\rfloor}$, so that $i\in \Gamma_{n_0+1,\lfloor p/8\rfloor}^p$.

\smallskip

Hence, for $p$ large enough, the set $\b_{n_0}$ is covered by at most $\#\Gamma_{n_0+1,\lfloor p/8\rfloor}^p $
sets of the form $\bp([(i-1)2^{-p},i2^{-p}])$, whose diameter is bounded above by
$C(\omega)(1+\log 2^p)2^{-p/4}$ thanks to Lemma \ref{Holder}. By the very definition of Hausdorff measures, we get
$$m^h(\b_{n_0})\leq \liminf_{p\to\infty} \Big(\#\Gamma_{n_0+1,\lfloor p/8\rfloor}^p \times h\Big(C(\omega)(1+\log 2^p)2^{-p/4}\Big)\Big)=0,$$
by \eqref{upper-Haus2}. This completes the proof of the first assertion of the proposition.

To get the second assertion, let $A$ be a Borel subset of $\bm_\infty$, and 
set
$$A':=\Bigg\{a\in A: \limsup_{r\downarrow 0}\frac{\mathrm{Vol}(B(a,r))}{h(r)} \geq K_2\Bigg\},$$
so that $m^h(A')=m^h(A)$ (by the first assertion) and $\mathrm{Vol}(A')\leq \mathrm{Vol}(A)$. 
Applying Lemma \ref{compa}(ii) to $A'$ gives the desired result.
\endproof

\section{Proof of the main result}
\label{proof-main}

In this section, we prove Theorem \ref{main}. We argue under the
probability measure $\N^{(1)}_0$. We start with a lemma.

\begin{lemma}
\label{nointersec}
$\N^{(1)}_0(\dd\omega)$ almost surely, for every $s\in(0,1)$, we
have
\begin{equation}
\label{no-single}
m^h(\bp([0,s])\cap\bp((s,1]))=\mathrm{Vol}(\bp([0,s])\cap\bp((s,1]))=0.
\end{equation}
\end{lemma}

\proof If $s\in(0,1)$,
a point $a$ of $\bm_\infty$ belongs to $\bp([0,s])\cap\bp((s,1])$ if and only if
there exist $r\in [0,s]$ and $r'\in(s,1]$ such that $a=\bp(r)=\bp(r')$, and in particular the
equivalence class of $r$, or of $r'$, for $\approx_{(\omega)}$ is not a singleton.
As a consequence of \cite[Theorem 3.4]{Invent}, it is almost surely true 
that, for any $r\in(0,1)$, the equivalence class 
of $r$ for $\approx_{(\omega)}$ is not a singleton only if either $p_{(\omega)}(r)$ is not a leaf
of the tree $\t_{(\omega)}$ (equivalently there exists $r'\in[0,1]\backslash \{r\}$
such that $p_{(\omega)}(r)=p_{(\omega)}(r')$), or if there exists $\ve>0$ such that
$\wh W_u\geq \wh W_r$ for every $u\in[r,(r+\ve)\wedge 1]$, or for
every $u\in[(r-\ve)\vee 0,r]$. The set of all $r\in[0,1]$ such that 
$p_{(\omega)}(r)$ is not a leaf
of the tree $\t_{(\omega)}$ has Lebesgue measure $0$. Similarly, 
for every fixed $r\in(0,1)$ and $\ve>0$, the property $\wh W_u\geq \wh W_r$ for every $u\in[r,(r+\ve)\wedge 1]$,
resp.~for every $u\in[(r-\ve)\vee 0,r]$), holds with zero probability under
$\N^{(1)}_0$, as an easy application of the properties of the Brownian snake (we omit the
details). We conclude that we have $\N^{(1)}_0$ a.s. for every $s\in(0,1)$,
$$\int_0^1 \dd r\,\mathbf{1}_{\{\bp(r)\in\bp([0,s])\cap\bp((s,1])\}} = 0,
$$
or equivalently 
$$\mathrm{Vol}(\bp([0,s])\cap\bp((s,1]))=0.$$
Proposition \ref{upper-Haus} then allows us to replace $\mathrm{Vol}$ by $m^h$ in the last display. \endproof

We denote Lebesgue measure on $[0,1]$ by $\mathscr{L}(\dd s)$. Together with Propositions \ref{lower-Haus} and \ref{upper-Haus},
the following lemma is the key ingredient
of the proof of Theorem \ref{main}. 

\begin{lemma}
\label{zero-one}
There exists a constant $\kappa\in[0,\infty]$ such that
we have $\N_0^{(1)}(\dd \omega)$ a.s., for $\mathscr{L}$ almost every $s\in(0,1)$,
$$\limsup_{\ve\to 0} \frac{m^h(\bp([s-\ve,s+\ve]))}{2\ve} =\kappa.$$
\end{lemma}

We postpone the proof of Lemma \ref{zero-one} and complete the proof
of Theorem \ref{main}. From now on, we fix $\omega$
such that the conclusions of Propositions \ref{lower-Haus} and \ref{upper-Haus}
and of Lemmas \ref{nointersec} and \ref{zero-one} hold.
For every Borel subset $A$ of $\bm_\infty$, we 
set $\mu(A)=m^h(A)$ to simplify notation. By \cite[Theorem 27]{Rog}, $\mu$ is a positive measure on the
Borel $\sigma$-field of $\bm_\infty$, and Proposition \ref{upper-Haus} entails that
$\mu$ is finite and absolutely continuous with respect to $\mathrm{Vol}$. More precisely,
Propositions \ref{lower-Haus} and \ref{upper-Haus} show that we can choose the
Radon-Nikodym derivative 
$$g:=\frac{\dd \mu}{\dd \mathrm{Vol}}$$
so that $\kappa_1\leq g\leq \kappa_2$. To prove Theorem \ref{main}, it will be enough to verify that $g=\kappa$, $\mathrm{Vol}$ a.e.,
where $\kappa$ is the constant in Lemma \ref{zero-one} --- note that this will automatically entail that $\kappa_1\leq \kappa\leq \kappa_2$.

For every $s\in [0,1]$, set $\Phi(s):=\mu(\bp([0,s]))$. Then, we have a.s. for every
$s\in[0,1]$,
$$
\Phi(s)=\int_{\bp([0,s])} g\,\dd \mathrm{Vol}
=\int_0^1 \mathbf{1}_{\{\bp(r)\in\bp([0,s])\}}\,g(\bp(r))\,\dd r
=\int_0^s g(\bp(r))\,\dd r,
$$
where the second equality is the definition of $\mathrm{Vol}$, and the third one holds because
$$\int_s^1 \mathbf{1}_{\{\bp(r)\in\bp([0,s])\}}\,\dd r\leq \int_0^1 \mathbf{1}_{\{\bp(r)\in\bp([0,s])\cap\bp((s,1])\}}\,\dd r=\mathrm{Vol}(\bp([0,s])\cap\bp((s,1]))=0$$
by \eqref{no-single} (we assumed $0<s<1$ in \eqref{no-single}, but the
cases $s=0$ and $s=1$ are trivial). 
From a classical result about differentiability of absolutely
continuous real functions, we infer that $\Phi$ is a.e.~differentiable, and
moreover, $\Phi'(s)=g(\bp(s))$, $\mathscr{L}(\dd s)$ a.e. Consequently,
we have
$$g(\bp(s))=\lim_{\ve\downarrow 0} \frac{\Phi(s+\ve)-\Phi(s-\ve)}{2\ve},$$
for $\mathscr{L}$ almost every $s\in(0,1)$. However,
$$\Phi(s+\ve)-\Phi(s-\ve)= \mu(\bp([0,s+\ve]))-\mu(\bp([0,s-\ve]))= \mu(\bp([s-\ve,s+\ve])),$$
using Lemma \ref{nointersec} in the last equality.

Then it follows from the conclusion of Lemma \ref{zero-one}  that we have a.s.
$g(\bp(s))=\kappa$ for $\mathscr{L}$ almost every $s\in(0,1)$.
Equivalently, $\int_0^1 \mathbf{1}_{\{g(\bp(s))=\kappa\}}\,\dd s=1$, 
which exactly means by the definition of $\mathrm{Vol}$ that
$$\int_{\bm_\infty} \mathbf{1}_{\{g(a)=\kappa\}}\,\mathrm{Vol}(\dd a) =1,$$
and thus $g=\kappa$, $\mathrm{Vol}$ a.e., which was the desired result. This completes the proof of Theorem \ref{main}. 

\medskip
\proof[Proof of Lemma \ref{zero-one}] 
In order to verify that the limsup  in Lemma \ref{zero-one} is constant, say for a fixed value of $s$, 
one would like to argue that this limsup
is measurable with respect to 
a trivial $\sigma$-field containing only events of probability zero or one. However, because 
of the rather involved construction of the Brownian sphere, it is not so easy to
find such a $\sigma$-field, and we will proceed in a slightly different manner. As in preceding
proofs, it will be convenient to replace $\N^{(1)}_0$ by the excursion measure $\N_0$, or rather by a
size-biased version of $\N_0$. More precisely, we introduce the $\sigma$-finite
measure $\N_0^\bullet$ defined on $\S_0\times \R_+$ by
$$\N_0^\bullet(\dd \omega\,\dd r)= \mathbf{1}_{\{0\leq r\leq \sigma(\omega)\}}\,\N_0(\dd\omega)\,\dd r.$$
For $(\omega,r)\in \S_0\times \R_+$, we use the notation $U(\omega,r)=r$. The
construction of Section \ref{Brown-sphere} allows us to associate a
compact metric space $\bm_\infty(\omega)=[0,\sigma(\omega)]\,/\!\approx_{(\omega)}$
with any $\omega\in\S_0$ (and we keep  the notation $\bp_{(\omega)}$ for the 
canonical projection from $[0,\sigma(\omega)]$ onto $\bm_\infty(\omega)$). As previously
we drop $\omega$ in the notation. By simple scaling arguments left to
the reader, the proof of Lemma \ref{zero-one} reduces to finding a 
constant $\kappa\in[0,\infty]$ such that
\begin{equation}
\label{zero-one-tech}
\N_0^\bullet\hbox{ a.e., } \limsup_{\ve\to 0} \frac{m^h(\bp([U-\ve,U+\ve]))}{2\ve} =\kappa.
\end{equation}
Note that using $\N_0^\bullet$ and $U$ in \eqref{zero-one-tech}
yields an asymptotic result valid for Lebesgue almost every $r$ in $[0,\sigma]$ under $\N_0$,
which is what we need to get Lemma \ref{zero-one}.

To derive \eqref{zero-one-tech}, we will rely on an appropriate zero-one law. %Let us argue under $\N_0^\bullet$. 
By the Bismut decomposition theorem of the Brownian excursion
(see e.g. \cite[Theorem XII.4.7]{RY}), the distribution of $\zeta_U$
under $\N_0^\bullet$ is Lebesgue measure on $\R_+$, and conditionally on $\zeta_U=x>0$,
the two processes $\zeta'_s:=\zeta_{(U-s)\vee 0}$
and $\zeta''_s:=\zeta_{(U+s)\wedge \sigma}$ are distributed as two independent linear
Brownian motions started from $x$ and stopped upon hitting $0$. Let us also define,
for every $s\geq 0$,
$$V'_s:= \wh W_{(U-s)\vee0} - \wh W_U\;,\quad V''_s:= \wh W_{(U+s)\wedge \sigma} - \wh W_U,$$
and set, for every $y>0$,
$$\tau'_y:=\inf\{s\geq 0: \zeta'_s=(\zeta_U- y)\vee 0\}\;,\quad \tau''_y:=\inf\{s\geq 0: \zeta''_s=(\zeta_U- y)\vee 0\}.$$
Finally, let $\g_y$ be the $\sigma$-field generated by the two pairs
$$\Big( \zeta'_s-\zeta_U, V'_s\Big)_{0\leq s\leq \tau'_y}\;,\ \Big( \zeta''_s-\zeta_U, V''_s\Big)_{0\leq s\leq \tau''_y}.$$
Fix $\chi \in(0,1)$ and write $\N^{\bullet,\chi}_0$ for the probability measure $\N^{\bullet}_0(\cdot\mid \chi <\zeta_U<1/\chi)$. We observe
that the tail $\sigma$-field
$$\g_{0+}:=\bigcap_{0<y<\chi} \g_y$$
is $\N^{\bullet,\chi}_0$-trivial, in the sense that it contains only events of $\N^{\bullet,\chi}_0$-probability zero or one. This is 
an easy application of properties of the Brownian snake and we only sketch the argument. First we set $\xi_t=W_U((\zeta_U-t)\vee 0)-\wh W_U$
for $0\leq t\leq \chi$, in such a way that $(\xi_t)_{0\leq t\leq \chi}$ is a linear Brownian motion started from $0$
under $\N^{\bullet,\chi}_0$. Then, in a way very similar to \cite[Lemma V.5]{Zurich}, we can construct a pair $(\n',\n'')$ of 
point measures on $\S_0\times [0,\chi]$, which (under $\N^{\bullet,\chi}_0$)  are independent Poisson measures
on $\S_0\times [0,\chi]$ with intensity $2\,\N_0(\dd\omega)\,\dd t$, and are also independent of $\xi$, in such a 
way that, for every $0<y<\chi$, $\g_y$ is the $\sigma$-field generated by $(\xi_t)_{0\leq t\leq y}$ and the restrictions
of $\n'$ and $\n''$ to $\S_0\times [0,y]$ (informally, while the process $\xi$ gives the variation of ``labels'' along the part of the ancestral 
line of $p_{(\omega)}(U)$ between heights $\zeta_U-\chi$ and $\zeta_U$ in the tree $\t_{(\omega)}$, the point measures $\n'$ and $\n''$ 
correspond to the subtrees branching off the left side and the right side of this part of the
ancestral line of $p_{(\omega)}(U)$). The triviality of $\g_{0+}$ immediately follows from
this description of the $\sigma$-fields $\g_y$. 

So to complete the proof of the lemma, we only need to verify that the limsup in \eqref{zero-one-tech}
is $\g_{0+}$-measurable. To  this end, we first introduce some notation. We argue under $\N^{\bullet,\chi}_0$ and fix $y\in(0,\chi)$. We note that we can uniquely define a snake trajectory $\omega^{(y)}\in\S_0$
with duration $\sigma(\omega^{(y)})=\tau'_y+\tau''_y$,
by the prescriptions
$$\zeta_s(\omega^{(y)}):=\left\{\begin{array}{ll}
\zeta'_{\tau'_y-s}-\zeta_U+y,\quad &\hbox{if }0\leq s\leq \tau'_y,\\
\noalign{\smallskip}
\zeta''_{s-\tau'_y}-\zeta_U+y,\quad&\hbox{if }\tau'_y\leq s\leq \tau'_y+\tau''_y,
\end{array}
\right.
$$
and 
$$\wh W_s(\omega^{(y)}):=\left\{\begin{array}{ll}
V'_{\tau'_y-s}-V'_{\tau'_y} &\hbox{if }0\leq s\leq \tau'_y,\\
\noalign{\smallskip}
V''_{s-\tau'_y}-V''_{\tau''_y},\quad&\hbox{if }\tau'_y\leq s\leq \tau'_y+\tau''_y.
\end{array}
\right.
$$
Note that $V'_{\tau'_y} = V''_{\tau''_y}$ by the snake property. Informally, the snake trajectory
$\omega^{(y)}$ codes the part of the snake trajectory $\omega$
corresponding to the excursion of $(\zeta_s)_{s\geq 0}$ above level $\zeta_U-y$ that straddles $U$.
It is immediate that $\omega^{(y)}$ is $\g_y$-measurable. With $\omega^{(y)}$,
we associate the metric space $\bm_\infty^{(y)}=\bm_\infty(\omega^{(y)})$ via the construction 
of Section \ref{Brown-sphere}, and we write $\bp^{(y)}$ for the canonical projection
from $[0,\sigma(\omega^{(y)})]$ onto $\bm_\infty^{(y)}$. We also let $D^{(y)}=D_{(\omega^{(y)})}$
be the metric on $\bm_\infty^{(y)}$. 

\smallskip
\noindent{\it Claim.} $\N^{\bullet,\chi}_0$ a.s., for every $\ve>0$ small enough, the metric space $\bp([U-\ve,U+\ve])$
equipped with the metric $D$ is isometric to the space $\bp^{(y)}([\tau'_y-\ve,\tau'_y+\ve])$ equipped 
with the metric $D^{(y)}$.

\smallskip
Assuming that the claim holds, we can replace $m^h(\bp([U-\ve,U+\ve]))$ by $m^h(\bp^{(y)}([\tau'_y-\ve,\tau'_y+\ve]))$ 
in \eqref{zero-one-tech} and
the fact that $\omega^{(y)}$ and $\tau'_y$ are both $\g_y$-measurable implies that the limsup in \eqref{zero-one-tech} is also
$\g_y$-measurable. Since this holds for every $y\in(0,\chi)$, the limsup is $\g_{0+}$-measurable 
and therefore is equal to a constant a.s.

It remains to prove our claim. We first choose a (random) real $z\in(0,y)$ 
such that $W_U(\zeta_U-z)<\wh W_U$ (this is possible since 
we already noted that $(W_U(\zeta_U-t)-\wh W_U)_{0\leq t\leq y}$ is a linear Brownian motion
started from $0$ under $\N^{\bullet,\chi}_0$). Notice that $W_U(\zeta_U-z)=\wh W_{U-\tau'_z}=\wh W_{U+\tau''_z}$
by the snake property. Set $\alpha=\frac{1}{2}(\wh W_U-W_U(\zeta_U-z))$ in such a way 
that 
\begin{equation}
\label{01-tech1}
\wh W_{U-\tau'_z}=\wh W_{U+\tau''_z} <\wh W_U-\alpha.
\end{equation}
Then we have 
$$\inf_{t\notin (U-\tau'_z,U+\tau''_z)} D(U,t)>\alpha$$
as a consequence of the so-called cactus bound (see e.g. formula (4) in \cite{Plane}), which implies that $D(U,t)$
is greater than or equal to $\wh W_U$ minus the minimum of labels along
the geodesic segment between $p_{(\omega)}(U)$ and $p_{(\omega)}(t)$ in
the tree $\t_{(\omega)}$ (note that this geodesic segment contains $p_{(\omega)}(U-\tau'_z)=p_{(\omega)}(U+\tau''_z)$
if $t\notin (U-\tau'_z,U+\tau''_z)$).
It follows that, for $\ve >0$ small enough,
we have also
\begin{equation}
\label{01-tech2}
\build{\inf_{s\in[U-\ve,U+\ve]}}_{t\notin (U-\tau'_z,U+\tau''_z)}^{} D(s,t)>\alpha.
\end{equation}
%By choosing $\alpha$ smaller if necessary, we can assume that
%\begin{equation}
%\label{01-tech2}
%\min_{s\in[U-\tau'_y,U]} \wh W_s - \wh W_U <-\alpha\;,\quad \min_{s\in[U,U+\tau''_y]} \wh W_s - \wh W_U <-\alpha.
%\end{equation}
On the other hand, if $\ve$ is small, we have
\begin{equation}
\label{01-tech3}
\sup_{s,t\in[U-\ve,U+\ve]} D(s,t) <\alpha/2
\end{equation}
which implies in particular
\begin{equation}
\label{01-tech4}
\sup_{s\in[U-\ve,U+\ve]} |\wh W_s-\wh W_U| <\alpha/2.
\end{equation}
Suppose now that $\ve$ is small enough so that the bounds of the preceding displays hold. When
applying formula \eqref{formula-D} to evaluate $D(s,t)$ for $s,t\in[U-\ve,U+\ve]$, we
may restrict our attention to choices of $s_i,t_i$ that belong to 
$[U-\tau'_z,U+\tau''_z]$ (indeed, if for instance $s_j\notin [U-\tau'_z,U+\tau''_z]$, we have
$\sum_{i=1}^j D^\circ(t_{i-1},s_i) \geq D(s,s_j)> \alpha$ by \eqref{01-tech2}, whereas
$D(s,t)<\alpha/2$ by \eqref{01-tech3}). Furthermore, if all reals $s_i,t_i$ belong to 
$[U-\tau'_z,U+\tau''_z]$, we may also assume that, for every $j\in\{1,\ldots,p\}$, the maximum appearing
 in formula \eqref{formula-D-0} for $D^\circ(t_{j-1},s_j)$ is 
$\min_{u\in[t_{i-1}\wedge s_j,t_{j-1}\vee s_j]}\wh W_u$. In fact, if this is not the case, this maximum 
is $\min_{u\in[t_{j-1}\vee s_j,t_{j-1}\wedge s_j]}\wh W_u$, which is smaller than $\wh W_U-\alpha$
by \eqref{01-tech1}, and then \eqref{formula-D-tech} and \eqref{01-tech4}
imply that $\sum_{i=1}^p D^\circ(t_{i-1},s_i) >\alpha$. Similar considerations apply to the evaluation
of $D^{(y)}(s,t)$:
%(note that the distribution of $(\omega^{(y)},\tau'_y)$ is $\N^\bullet_0(\cdot\mid \zeta_U=y)$): 
With the same
value of $z$, we obtain that, when $s,t\in[\tau'_y-\ve,\tau'_y+\ve]$ and $\ve$ is small, the infimum in the analog of
formula \eqref{formula-D} giving $D^{(y)}(s,t)$ can be restricted to $s_i,t_i\in [\tau'_y-\tau'_z,\tau'_y+\tau''_z]$, and to the case
where, for every $j$, the maximum in formula \eqref{formula-D-0} for $D^\circ_{(\omega^{(y)})}(t_{j-1},s_j)$
is equal to $\min_{u\in[t_{i-1}\wedge s_j,t_{j-1}\vee s_j]}\wh W_u(\omega^{(y)})$. 
 Finally, we also observe that,
for every $s,t\in [U-\tau'_y,U+\tau''_y]$, the property $p_{(\omega)}(s)=p_{(\omega)}(t)$
holds if and only if $p_{(\omega^{(y)})}(\tau'_y+(s-U))=p_{(\omega^{(y)})}(\tau'_y+(t-U))$.

It now follows from the preceding considerations that, for $\ve$ small enough, for
every $s,t\in[U-\ve,U+\ve]$, we have
$$D(s,t)=D^{(y)}(\tau'_y+(s-U),\tau'_y+(t-U)).$$
Indeed, each term in the infimum appearing in formula \eqref{formula-D}  for $D(s,t)$  with a choice of
$s_i,t_i\in [U-\tau'_z,U+\tau''_z]$, corresponds to an analogous term in the formula
for $D^{(y)}(\tau'_y+(s-U),\tau'_y+(t-U))$, with the choice $s'_i=\tau'_y + (s_i-U), t'_i=\tau'_y+(t_i-U)$, 
and these two terms are immediately seen to be equal. 
We conclude that the mapping $s\mapsto \tau'_y+(s-U)$ induces an isometry from 
$\bp([U-\ve,U+\ve])$ onto $\bp^{(y)}([\tau'_y-\ve,\tau'_y+\ve])$. This completes the
proof of the claim and of Lemma \ref{zero-one}. \endproof

 \end{document}